\documentclass[12pt]{article}
\usepackage[dvips]{epsfig}
\usepackage{graphicx}
\usepackage{amsthm,amsmath,amssymb}
\usepackage{latexsym}
\usepackage{graphics}
\usepackage{url}
\usepackage{placeins}
\usepackage{hyperref}
\usepackage{color}
\usepackage{bm}
\usepackage{comment}
\usepackage{statex}

\definecolor{blue}{rgb}{0,0.2,1.0}
\definecolor{red}{rgb}{0.9,0,0}
\definecolor{green}{rgb}{0,0.50,0.10}
\definecolor{violet}{rgb}{0.5804,0.0000,0.8275}

\def\@themcountersep{}

\newtheorem{THEO}{Theorem}[section]
\newtheorem{ALGo}[THEO]{Algorithm}
\newtheorem{CONJ}[THEO]{Conjecture}
\newtheorem{COND}[THEO]{Condition}
\newtheorem{CORO}[THEO]{Corollary}
\newtheorem{DEFI}[THEO]{Definition}
\newtheorem{EXAMP}[THEO]{Example}
\newtheorem{FACT}[THEO]{Fact}
\newtheorem{HYPO}[THEO]{Hypothesis}
\newtheorem{LEMM}[THEO]{Lemma}
\newtheorem{PROB}[THEO]{Problem}
\newtheorem{PROP}[THEO]{Proposition}
\newtheorem{REMA}[THEO]{Remark}
\newcommand{\theo}{\begin{THEO}}
\newcommand{\algo}{\begin{ALGo} \rm}
\newcommand{\cond}{\begin{COND}}
\newcommand{\conj}{\begin{CONJ}}
\newcommand{\coro}{\begin{CORO}}
\newcommand{\defi}{\begin{DEFI} \rm}
\newcommand{\examp}{\begin{EXAMP} \rm}
\newcommand{\fact}{\begin{FACT}}
\newcommand{\hypo}{\begin{HYPO} \rm}
\newcommand{\lemm}{\begin{LEMM}}
\newcommand{\prob}{\begin{PROB} \rm}
\newcommand{\prop}{\begin{PROP}}
\newcommand{\rema}{\begin{REMA} \rm}
\newcommand{\etheo}{\end{THEO}}
\newcommand{\ealgo}{\end{ALGo}}
\newcommand{\econd}{\end{COND}}
\newcommand{\econj}{\end{CONJ}}
\newcommand{\ecoro}{\end{CORO}}
\newcommand{\edefi}{\end{DEFI}}
\newcommand{\eexamp}{\end{EXAMP}}
\newcommand{\efact}{\end{FACT}}
\newcommand{\ehypo}{\end{HYPO}}
\newcommand{\elemm}{\end{LEMM}}
\newcommand{\eprob}{\end{PROB}}
\newcommand{\eprop}{\end{PROP}}
\newcommand{\erema}{\end{REMA}}

\def\0{\mbox{\bf 0}}
\def\1{\mbox{\bf 1}}
\def\2{\mbox{\bf 2}}
\def\3{\mbox{\bf 3}}
\def\4{\mbox{\bf 4}}
\def\5{\mbox{\bf 5}}
\def\6{\mbox{\bf 6}}
\def\7{\mbox{\bf 7}}
\def\8{\mbox{\bf 8}}
\def\9{\mbox{\bf 9}}
\def\a{\mbox{\boldmath $a$}}
\def\b{\mbox{\boldmath $b$}}
\def\cc{\mbox{\boldmath $c$}}
\def\d{\mbox{\boldmath $d$}}
\def\e{\mbox{\boldmath $e$}}
\def\f{\mbox{\boldmath $f$}}

\def\s{\mbox{\boldmath $s$}}

\def\u{\mbox{\boldmath $u$}}

\def\x{\mbox{\boldmath $x$}}

\def\A{\mbox{\boldmath $A$}}
\def\B{\mbox{\boldmath $B$}}

\def\H{\mbox{\boldmath $H$}}
\def\I{\mbox{\boldmath $I$}}

\def\O{\mbox{\boldmath $O$}}

\def\Q{\mbox{\boldmath $Q$}}

\def\X{\mbox{\boldmath $X$}}
\def\Y{\mbox{\boldmath $Y$}}

\def\FC{\mbox{$\cal F$}}

\def\SC{\mbox{$\cal S$}}

\def\inprod#1#2{\langle#1, \, #2\rangle}

\def\Real{\mbox{$\mathbb{R}$}}

\def\s0{\mbox{\scriptsize \boldmath $0$}}

\def\sSC{\mbox{\scriptsize $\SC$}}

\def\Real{\mathbb{R}}
\def\coneC{\mathbb{C}}
\def\coneF{\mathbb{F}}
\def\coneK{\mathbb{K}}
\def\coneJ{\mathbb{J}}

\def\spaceV{\mathbb{V}}
\def\SymMat{\mathbb{S}}

\begin{document}

\title{Further Development in Convex Conic Reformulation of Geometric 
Nonconvex Conic Optimization Problems}

\author{
\normalsize
Naohiko Arima\thanks{
	\tt nao$\_$arima@me.com.}, \and \normalsize
Sunyoung Kim\thanks{Department of Mathematics, Ewha W. University, 52 Ewhayeodae-gil, Sudaemoon-gu, Seoul 03760, Korea 
			({\tt skim@ewha.ac.kr}). 
			 The research was supported  by   NRF 2021-R1A2C1003810.}, \and \normalsize
Masakazu Kojima\thanks{Department of Industrial and Systems Engineering,
	Chuo University, Tokyo 192-0393, Japan ({\tt kojima@is.titech.ac.jp}).
 	}
}


\maketitle 


\begin{abstract}
\noindent
A geometric nonconvex conic optimization problem (COP) 
 was recently proposed by Kim, Kojima and Toh (SIOPT 30:1251-1273, 2020) as
a unified framework for convex conic reformulation of a class of 
quadratic optimization problems 
and polynomial optimization problems.   
The nonconvex COP minimizes a linear function over the intersection of a nonconvex cone $\coneK$, 
a convex subcone $\coneJ$ of the convex hull co$\coneK$ of $\coneK$, 
and an affine hyperplane with a normal vector $\H$. Under the assumption 
co$(\coneK\cap\coneJ) = \coneJ$,
the original nonconvex COP in their paper was shown to be equivalently
formulated as a convex conic program by replacing the constraint set with the intersection of $\coneJ$  
and the affine hyperplane. 
This paper further studies 
the key assumption co$(\coneK \cap \coneJ) = \coneJ$ 
in their framework  and provides three sets of necessary-sufficient conditions 
for the assumption. 
Based on the conditions, 
we propose a new wide class of 
quadratically constrained quadratic programs with multiple nonconvex equality and inequality 
constraints, which can be  solved exactly by their semidefinite relaxation. 
\end{abstract}

\vspace{0.5cm}


\vspace{0.5cm}


\section{Introduction}

Convex conic 
reformulation of 
a geometric nonconvex conic optimization problem (COP) was studied by Kim, Kojima and Toh  
in  \cite{KIM2020}  as a unified framework for 
completely positive programming  reformulation of 
a wide class of nonconvex quadratic optimization problems (QOPs).
This class includes a  wide range of QOPs, such as QOPs 
over the standard simplex \cite{BOMZE2000}, maximum stable set problems \cite{DEKLERK2002}, 
graph partitioning problems \cite{POVH2007} and quadratic assignment problems 
\cite{DEKLERK2002},  
and its extension to polynomial optimization problems (POPs) \cite{ARIMA2013,ARIMA2014b,PENA2015}.  
The class of QOPs also covers Burer's class of QOPs with linear equality 
and complementarity constraints in nonnegative and binary variables \cite{BURER2009}. 
Their geometric nonconvex COP denoted by COP($\coneK\cap\coneJ$) 
below is quite simple, but it is powerful enough to capture 
the basic essentials to investigate 
convex conic reformulation
of general nonconvex COPs. 
In fact, it minimizes a linear function in a finite dimensional vector space $\spaceV$ 
over the intersection of three geometrically represented sets, 
a nonconvex cone $\coneK$, a convex subcone $\coneJ$ of the convex hull co$\coneK$ 
of $\coneK$, and an affine hyperplane with a normal vector~$\H$. 
The framework was also used  in the recent paper \cite{GABL2023} which proposed a large class of 
quadratically constrained 
quadratic programs with stochastic data. See also \cite{BOMZE2023}.

Throughout the paper, $\spaceV$ denotes a finite dimensional vector space endowed 
with the inner product $\inprod{\A}{\B}$ for every pair of $\A$ and $\B$ 
in $\spaceV$, and $\O \not= \H \in \spaceV$ and $\Q \in \spaceV$ are arbitrarily fixed, unless 
specified. 
For every cone
$\coneC \subseteq \spaceV$, let COP($\coneC$) denote the conic optimization 
problem (COP) of the form 
\begin{eqnarray}
\zeta_p(\coneC) = \inf\left\{\inprod{\Q}{\X} : \X \in \coneC, \ \inprod{\H}{\X} = 1 \right\}. 
\label{eq:COPC}
\end{eqnarray}
Here we say that 
$\coneC  \subseteq  \spaceV$ is a {\em cone}, which is not necessarily convex nor closed,
if $\lambda \A \in \coneC$ for every $\A \in \coneC$ and $\lambda \geq 0$. 
If COP($\coneC$) is infeasible, we let $\zeta_p(\coneC) = \infty$. 
By replacing $\coneC$ by its convex hull co$\coneC$, we obtain {\em a convex relaxation} 
COP(co$\coneC$) of the problem. 
When the original nonconvex COP($\coneC$)  and the relaxed convex COP(co$\coneC$) 
share a common optimal value, {\it i.e.}, $\zeta_p(\coneC) = \zeta_p(\mbox{co}\coneC)$, 
the convex COP is called {\em a convex conic reformulation} 
of the nonconvex COP.  

A nonconvex COP introduced in \cite{KIM2020} is described as 
\begin{eqnarray*}
& & \mbox{COP($\coneK\cap\coneJ$):}\quad \zeta_p(\coneK\cap\coneJ) = \inf\left\{\inprod{\Q}{\X} : \X \in \coneK\cap\coneJ, \ \inprod{\H}{\X} = 1 \right\}
\end{eqnarray*}
under the following conditions: 
\vspace{2mm}\\
\noindent
{\bf (A) } $\coneK$ is a nonempty cone in $\spaceV$ and 
$-\infty < \zeta_p(\coneK\cap\coneJ) < \infty$ ({\it i.e.}, COP($\coneK\cap\coneJ$) is feasible and has 
a finite optimal value). 
\vspace{1mm}\\
{\bf (B)} $\coneJ$ is a 
convex cone contained in co$\coneK$ 
such that co$(\coneK\cap\coneJ) = \coneJ$. 
\vspace{2mm}\\
Among the theoretical results established in \cite{KIM2020}, we mention the 
following equivalence result (see \cite[Theorem 3.1]{KIM2020}, \cite[Theorem 5.1]{ARIMA2022} for more details). 
\theo \label{theorem:mainKKT}  
Assume that Conditions (A) 
and (B) are satisfied.  Then,  
\begin{eqnarray}
\left. 
\begin{array}{c}
-\infty < \zeta_p(\coneJ) < \infty,\\
\Updownarrow \\
-\infty < \zeta_p(\coneK\cap\coneJ) = \zeta_p(\coneJ)  < \infty \\ 
\mbox{({\it i.e.}, COP($\coneJ$) is a convex reformulation of COP($\coneK\cap\coneJ$))}.
\end{array}
\right\} 
 \label{eq:equivCond}
\end{eqnarray}
\etheo

First, we briefly discuss some remaining issues, which were not throughly studied 
in \cite{KIM2020}: 
\begin{description}
\item{(a) } 
The feasibility preserving property \cite{YILDIRIM2023} (Section 2.1). 
\item{(b) } Strong duality of the reformulated convex COP (Section 2.2).
\item{(c) } The existence of a common optimal solution of a nonconvex COP and the 
reformulated convex COP (Section 2.3). 
\end{description}

If $\coneJ$ is a face of co$\coneK$, then co$(\coneK\cap\coneJ) = \coneJ$ 
(Condition (B)) is satisfied. This case was 
throughly studied and played an essential role in convex conic 
reformulation of QOPs and POPs in \cite{KIM2020}. Another case mentioned in 
\cite[Lemma 2.1 (iv)]{KIM2020} for Condition (B) is: 
$\coneJ$ is the convex hull of the union of 
(possibly infinitely many) faces of co$\coneK$. But neither its implication 
nor its application was discussed there. In Section 3 of this paper, we 
introduce the family $\widehat{\FC}(\coneK)$ of all $\coneJ$ which 
satisfy co$(\coneK\cap\coneJ) = \coneJ$, and study its fundamental properties.
In particular, we establish the following characterization of $\widehat{\FC}(\coneK)$. 
\theo \label{theorem:NS} \mbox{\ }
\begin{description}
\item{(i) } 
$\coneJ \in \widehat{\FC}(\coneK)$ iff $\coneJ = \mbox{co}\coneK'$ 
for some cone $\coneK' \subseteq \coneK$. 
\item{(ii) } 
$\coneJ \in \widehat{\FC}(\coneK)$ iff $\coneJ = \mbox{co}\big(\bigcup \FC \big)$ for some 
$\FC \subseteq \widehat{\FC}(\coneK)$.
\item{(iii) } Assume that $\coneK$ and $\coneJ$ are closed 
and $\H \in \mbox{int}(\coneK^*)$, 
where $\mbox{int}(\coneK^*)$ denotes  
the interior of $\coneK^*$. 
Then $\coneJ \in \widehat{\FC}(\coneK)$ iff 
$\zeta_p(\coneK \cap \coneJ) = \zeta_p(\coneJ)$ for every $\Q \in \spaceV$. 
\end{description}
\etheo
\noindent
A proof is given in Section~\ref{Subsection:proofNS}. 
Based on assertion (ii),   
we discuss a decomposition of the convex reformulation 
COP$(\mbox{co}(\bigcup\FC))$ of 
COP$(\coneK\cap\big(\mbox{co}(\bigcup\FC)\big))$ with 
$\FC \subseteq \widehat{\FC}(\coneK)$ 
into the convex reformulations COP($\coneF$) of COP($\coneK\cap\coneF$) $(\coneF \in \FC)$ (Theorem~\ref{theorem:equivalence}).

In Section~\ref{Section4}, we focus on the 
case where $\coneK$ is represented as $\coneK = \{\x\x^T : \x \in \Real^n\} \subseteq \SymMat^n$ (the space of $n \times n$ symmetric matrices) 
and $\coneJ$ by multiple inequalities such that 
$\coneJ = \coneJ_- \equiv 
\{ \X \in \mbox{co}\coneK : \inprod{\B_k}{\X} \leq 0 \ (1 \leq k \leq m)\}$ 
for some $\B_k \in \SymMat^n$ $(1 \leq k \leq m)$. 
In this case, co$\coneK$ forms the cone $\SymMat^n_+$ of positive 
semidefinite matrices in $\SymMat^n$, 
and COP$(\coneJ_-)$ serves as a semidefinite programming (SDP) relaxation of COP$(\coneK\cap\coneJ_-)$, 
which can be regarded as a quadratically constrained quadratic program (QCQP) with 
nonconvex inequality constraints $\x^T\B_k\x \leq 0$ $(1 \leq k \leq m)$ and 
an equality constraint $\x^T\H\x = 1$. 
By using Theorem~\ref{theorem:NS} (iii) and \cite[Lemma 2.2]{YE2003},  
we establish that 
$\coneJ_-\equiv \{ \X \in \mbox{co}\coneK : \inprod{\B_k}{\X} \leq 0 \ (1 \leq k \leq m)\}  \in \widehat{\FC}(\coneK)$ if condition 
\begin{eqnarray}
\coneJ_0(\B_k) \equiv \{\X \in \SymMat^n_+ : \inprod{\B_{k}}{\X} = 0\} 
\subseteq \coneJ_-\ (1 \leq k \leq m) \label{eq:condition0}
\end{eqnarray} 
is satisfied.  
This result leads to 
a wide class of QCQPs with multiple nonconvex constraints, COP$(\coneK\cap\coneJ_-)$ 
with $\B_k \in \SymMat^n$ $(1 \leq k \leq m)$,  
which can be solved by their SDP relaxation COP$(\coneJ_-)$. 
See Figure 1 for a geometrical image of condition~\eqref{eq:condition0}. We know that if $\X$ is a common 
optimal solution of COP$(\coneK\cap\coneJ_-)$ and its SDP relaxation COP$(\coneJ_-)$, then rank$\X = 1$. 
In  case (a), every extreme ray of $\coneJ_-$ 
on which COP$(\coneJ_-)$ can attain an optimal solution $\X$ lies on the boundary of 
$\SymMat^n_+$ whose extreme rays are known to be generated by 
rank-1 matrices. In case (b), $\coneJ_-$ includes an extreme ray not included on 
the boundary of $\SymMat^n_+$; 
hence such an extreme ray may include  
a matrix with rank greater than $1$.  
Thus, condition~\eqref{eq:condition0} is quite natural to ensure the equivalence 
of COP$(\coneJ_-)$ and its SDP relaxation COP$(\coneK\cap\coneJ_-)$ for any $\Q \in \SymMat^n$. 
We also see that $\coneJ_- = \mbox{co}\coneK'$ for some $\coneK' \subseteq \coneK$ in case (a) but 
$\coneJ_- \not= \mbox{co}\coneK'$ for any $\coneK' \subseteq \coneK$ in case (b). 
Therefore, by Theorem~\ref{theorem:NS} (i),  $\coneJ_- \in \widehat{\FC}(\coneK)$ in case (a) but $\coneJ_- \not\in \widehat{\FC}(\coneK)$ in case (b). 
We note that if $n \geq 3$  then the intersection of 
the boundary of $\SymMat^n_+$ with $\coneJ_0(\B)$ ($\B \in \SymMat^n$)  
generally includes  matrices with rank in $\{0,1,\ldots,n-1\}$, presenting a more complicated situation.
The aforementioned explanation represents the case when $n=2$.

\begin{figure}[t!]  
\begin{center}
\includegraphics[width=0.40\textwidth]{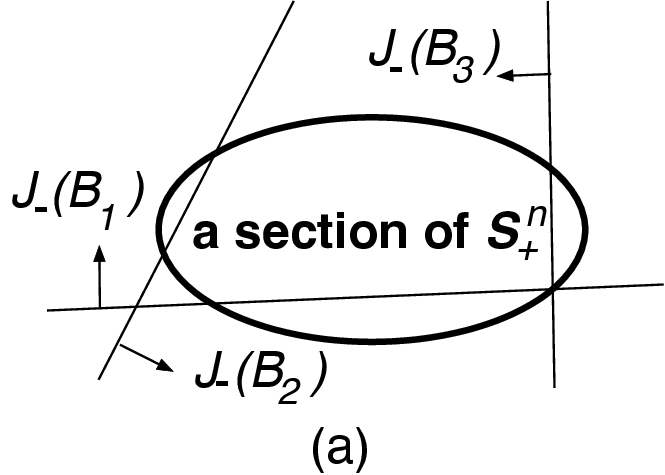}
\mbox{\ } \hspace{1cm}
\includegraphics[width=0.40\textwidth]{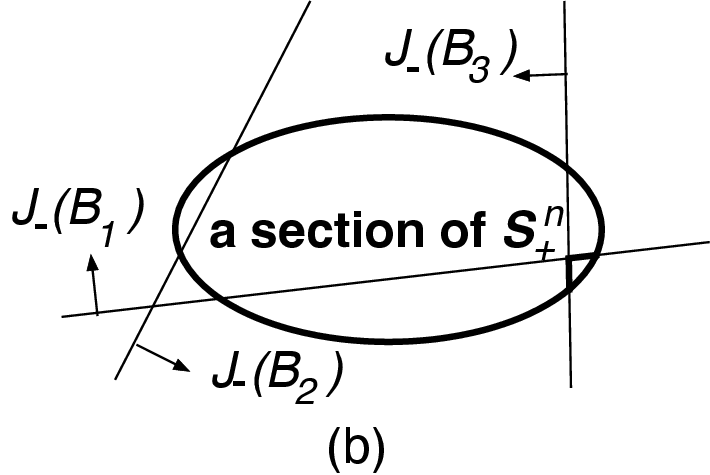}
\end{center}

\caption{Geometrical illustrations for condition~\eqref{eq:condition0} with $m=3$. 
(a) illustrates a 
case where~\eqref{eq:condition0} is satisfied, and (b) a case 
where~\eqref{eq:condition0}
is not satisfied. Here $\coneJ_-(\B_k) = \{ \X \in \SymMat^n_+ : 
\inprod{\B_k}{\X} \leq 0 \}$ $(1 \leq k \leq m)$ and 
$\coneJ_- = \cap_{k=1}^3\coneJ_-(\B_k)$.   
}
\end{figure}

\subsection{Contribution of the paper and related literature}

Our primary contribution lies in presenting 
a new and wide class of QCQPs with {\em multiple} nonconvex constraints 
that can be solved exactly by their SDP relaxation.  
The equivalence of QCQPs and their SDP relaxation was  extensively studied  
in the literature, 
including \cite{AZUMA2022,JEYAKUMAR2014,KIM2003,POLYAK1998,SOJOUDI2014,STERN1995,WANG2022,YE2003,ZHANG2000}. 
The classes of QCQPs studied there can be classified 
into two categories. The first category requires particular sign patterns 
of the data matrices involved 
in QCQPs \cite{AZUMA2022,KIM2003,SOJOUDI2014,ZHANG2000}. The studies on the second 
category have  focused on  cases where the number of nonconvex constraints is limited at most two 
\cite{YE2003},  and some additional assumption on the quadratic functions involved 
in the constraint \cite{JEYAKUMAR2014,POLYAK1998,STERN1995,WANG2022,YE2003}. 
Trust-region subproblems of nonlinear programs
 have been  frequently studied in relation to the  second category. 
Our class of QCQPs is also related to the second category, but 
represents a significantly wider class than them, in the sense that QCQPs can involve 
any finite number of nonconvex 
constraints. Condition~\eqref{eq:condition0} 
is both quite general and 
intuitively natural,  ensuring the equivalence of QCQPs and their SDP relaxation 
(see Figure 1).

The second contribution is that we present more precise 
characterizations of Condition (B), which plays a 
central role in the theory of convex conic reformulation of geometric COPs, 
as shown in Theorem~\ref{theorem:NS}. 
This contribution, together with (a), (b) and (c) mentioned above,  makes the theory 
solid and enhances its applicability   to a wide range of problems. 
In particular, the discovery of the new class of QCQPs 
in the first contribution is based on them.

\subsection{Outline of the paper}

In Section 2, we discuss the three issues (a), (b) and (c) 
in detail. In Section~3, we present some fundamental properties on 
$\widehat{\FC}(\coneK)$, and prove Theorem~\ref{theorem:NS}. 
We also discuss decompositions of a nonconvex COP and its convex conic reformulation 
based on Theorem~\ref{theorem:NS} (ii). 
In Section 4, we present the class of QCQPs with multiple 
nonconvex inequality and equality constraints that can be reformulated as their SDP relaxation. Some representative 
 QCQP examples in this class are presented. 
We conclude in Section 5 with remarks on the features of the geometric nonconvex COP$(\coneK \cap \coneJ)$ 
as a unified framework, 
and its possible application to the completely positive cone.

\subsection{Remarks on notation and symbols}

We note that a cone $\coneC$ is {\em convex} if $\X =\sum_{k=1}^m \X^k \in \coneC$ whenever 
$\X^k \in \coneC$  $(k=1,...,m)$.  
The convex hull co$\coneC$ of a cone $\coneC$ is represented as co$\coneC =\{\sum_{k=1}^m \X^k : \X^k \in \coneC \ (k=1,\ldots,m) \ \mbox{for some } m \}$. 
The dual $\coneC^* = \{ \Y  \in  \spaceV  : \inprod{\X}{\Y}  \geq  0 \ 
\mbox{for every } \X \in \coneC \}$
of a cone $\coneC$ always forms a closed convex cone in $\spaceV$. 
Let $\coneC$ be a convex cone in $\spaceV$. 
A convex cone $\coneF \subset \coneC$ is 
a {\em face} of $\coneC$ if $\X^k \in \coneF$ $(1 \leq k \leq m)$ 
whenever $\X =\sum_{k=1}^m \X^k \in \coneF$ 
and $\X^k \in \coneC$ $(1 \leq k \leq m)$.
An {\em extreme ray} of $\coneC$ 
is a face which spans a $1$-dimensional linear subspace of~$\spaceV$.


\section{Some remaining issues}

\label{Section2}

In this section, we discuss issues (a) and (b) and (c) mentioned in Section 1. 

\subsection{The feasibility preserving property - (a)} 

The property that the nonconvex problem is feasible iff its 
convex relaxation is feasible is called feasibility preserving. 
This property was fully studied in \cite{YILDIRIM2023} for convex relaxation of nonconvex QOPs. 
We present a simple proof that our geometric framework is 
feasibility preserving. 

\lemm \label{lemma:feasiblity}
Let $\coneC$ be a nonempty cone in $\spaceV$. 
COP$(\coneC)$ is feasible (infeasible, respectively) iff 
COP$(\mbox{co}\coneC)$ is feasible (infeasible, respectively). 
\elemm
\proof{It suffices to show that if COP(co$\coneC$) is feasible, then COP($\coneC$) is feasible. 
Let $\overline{\X}$ be a feasible solution of COP(co$\coneC$). Then, there exist 
$\X_k \in \coneC$ $(1 \leq k \leq m)$ such that $\overline{\X} = \sum_{k=1}^m\X_k$. Since 
$1 = \inprod{\H}{\overline{\X}} = \inprod{\H}{\sum_{k=1}^m\X_k}$,  
$\inprod{\H}{\X_k} > 0$ for some $k$. Let $\widehat{\X} = \X_k / \inprod{\H}{\X_k}$.
Then $\widehat{\X} \in \coneC$ and $\inprod{\H}{\widehat{\X}} = 1$. Hence 
$\widehat{\X}$ is a feasible solution of COP($\coneC$). 
\qed
}

By Lemma~\ref{lemma:feasiblity}, we can weaken Condition (A). 
\begin{description}
\item[(A)' ] $\coneK$ is a nonempty (not necessarily convex) cone in $\spaceV$.
\end{description}
\coro \label{corollary:feasiblity}
Assume Conditions (A)' and (B). Then~\eqref{eq:equivCond} holds.
\ecoro
\proof{If $-\infty < \zeta_p(\coneJ) < \infty$, then 
$-\infty <  \zeta_p(\coneJ) \leq \zeta_p(\coneK\cap\coneJ) < \infty$ holds by Lemma~\ref{lemma:feasiblity}. 
Hence Condition (A) is satisfied. Therefore,~\eqref{eq:equivCond} holds 
by Theorem~\ref{theorem:mainKKT}.
\qed
} 

By the corollary and the lemma above, we know under Conditions (A)' and (B) that if 
COP$(\coneJ)$ attains the finite optimal value $\zeta_p(\coneJ)$, then 
$\zeta_p(\coneK\cap\coneJ) = \zeta_p(\coneJ)$, and that if COP$(\coneJ)$ is 
infeasible, then so is COP$(\coneK\cap\coneJ)$.   

\subsection{Strong duality of COP($\coneJ$) - (b)
}

As a straightforward application of \cite[Theorem 2.1]{KIM2022} 
to the primal-dual pair COP($\coneJ$) and DCOP($\coneJ$) with $\H \in \coneJ^*$, 
we obtain Theorem~\ref{theorem:duality} below. 
Here the dual of COP$(\coneJ)$ is given by 
\begin{eqnarray}
\mbox{DCOP$(\coneJ)$:} \ \zeta_d(\coneJ) = \sup\{t : \Q - \H t \in \coneJ^* \}. 
\label{eq:dualCOP0}  
\end{eqnarray}

\theo \label{theorem:duality}
Assume that $\coneJ$  is a closed convex cone in $\spaceV$, $\H \in \coneJ^*$ and that 
$-\infty < \zeta_p(\coneJ) < \infty$ or  $-\infty < \zeta_d(\coneJ) < \infty$.  
Then the following assertions  hold. 
\begin{description}
\item{(i)} 
$-\infty < \zeta_p(\coneJ)=\zeta_d(\coneJ) < \infty$ holds. 
\item{(ii)}  
Dual DCOP($\coneJ$) has an optimal solution. 
\item{(iii)} 
The set of optimal solutions of primal COP($\coneJ$) is nonempty 
and bounded  iff $\Q - \H t \in \mbox{int}\coneJ^*$ for some $t \in \Real$. 
\item{(iv)} 
The set of optimal  solutions of 
primal COP($\coneJ$) is nonempty and unbounded if 
$\coneJ$ is not pointed and 
$\Q - \H t \in \mbox{relint}\coneJ^*$ for some $t \in \Real$, where $\mbox{relint}\coneJ^*$ 
denotes the relative interior of $\coneJ^*$ with respect to the linear subspace of $\spaceV$ 
spanned by $\coneJ^*$.
\end{description}
\etheo
We note that the condition ``$\coneJ$ is closed'' is natural when the 
existence of an optimal solution of COP($\coneJ$) is discussed, and the condition $\H \in \coneJ^*$ 
holds naturally 
when  QCQPs and POPs are converted into COP($\coneK\cap\coneJ$). 
See \cite[Sections 3.2, 4 and 5]{KIM2020}.

\subsection{
Existence of a common optimal solution of COP($\coneJ$) and COP($\coneK\cap\coneJ$) - (c)}

If $\X$ is an optimal solution of COP($\coneK\cap\coneJ$), then 
the identity $\zeta_p(\coneK\cap\coneJ) = \zeta_p(\coneJ)$ ensures that 
$\X$ is also an optimal solution of COP($\coneJ$)  since $\coneK\cap\coneJ \subseteq \coneJ$. In general, 
however, a nonconvex (or even convex) COP may have no optimal solution even when 
it has a finite optimal value. See, for example, \cite[Section 4]{KIM2023}. 
The following theorem provides a sufficient condition for a common optimal solution of 
COP($\coneK\cap\coneJ$) and COP($\coneJ$).  
\theo \label{theorem:commonSol} 
Assume that 
$\coneK$ is a closed cone in $\spaceV$, 
$\coneJ \subseteq \mbox{co}\coneK$ is a closed convex cone satisfying 
co$(\coneK\cap\coneJ)=\coneJ$, COP$(\coneK\cap\coneJ)$ is feasible, 
$\H \in \coneJ^*$, and that $\Q - \H t \in \mbox{int}\coneJ^*$ (the interior of $\coneJ^*$) 
for some $t \in \Real$. Then, 
\begin{eqnarray}
\left. 
\begin{array}{l}
-\infty < \zeta_d(\coneJ) = \zeta_p(\coneJ) = \zeta_p(\coneK\cap\coneJ) < \infty, \\[3pt]
\mbox{COP$(\coneK\cap\coneJ)$ and COP$(\coneJ)$ have a common optimal solution,}\\[3pt]
\mbox{DCOP$(\coneJ)$ has an optimal solution.}
\end{array}
\right\} \label{eq:commonSol}
\end{eqnarray}
\etheo 
The strict feasibility of DCOP$(\coneJ)$ ({\it i.e.}, Slater's constraint qualification 
$\Q - \H t \in \mbox{int}\coneJ^*$ for some $t$) is assumed here 
for the existence of solutions of COP$(\coneK\cap\coneJ)$ and COP$(\coneJ)$. 
It should be emphasized that none of the finite optimal values for COP
$(\coneK\cap\coneJ)$, 
COP$(\coneJ)$ and DCOP$(\coneJ)$ is assumed in advance.

\noindent
{\em Proof of Theorem~\ref{theorem:commonSol}}. 
We prove that COP($\coneJ$) and COP($\coneK\cap\coneJ$) have a common optimal solution $\X^*$. 
Choose a feasible solution $\overline{\X} \in \coneK\cap\coneJ$ of COP($\coneK\cap\coneJ$). 
We consider the level sets of COP($\coneK\cap\coneJ$) and COP($\coneJ$) given by 
\begin{eqnarray*}
S(\coneK\cap\coneJ) & \equiv & \left\{ \X \in \coneK\cap\coneJ : \inprod{\H}{\X} = 1, \ 
\inprod{\Q}{\X} \leq \inprod{\Q}{\overline{\X}} \right\}, \\
S(\coneJ) & \equiv & \left\{ \X \in \coneJ : \inprod{\H}{\X} = 1, \ 
\inprod{\Q}{\X} \leq \inprod{\Q}{\overline{\X}} \right\}. 
\end{eqnarray*}
Then $\overline{\X} \in S(\coneK\cap\coneJ) \subseteq S(\coneJ)$. Since $\coneK$ and 
$\coneJ$ are closed, both $S(\coneK\cap\coneJ)$ and $S(\coneJ)$ are closed. 
We will show that $S(\coneJ)$ is bounded.
Assume on the contrary that $S(\coneJ)$ 
is unbounded. Then, there exists a sequence 
$\{\X_k \in S(\coneJ)\}$ such that $\parallel \X_k\parallel \rightarrow \infty$ as $k \rightarrow \infty$. We may assume without loss of generality 
that $\X_k/\parallel \X_k \parallel \in \coneJ$ 
converges to $\Delta\X  \in \coneJ$ such that 
\begin{eqnarray*}
\Delta\X \in \coneJ, \ \parallel \Delta\X\parallel  = 1, \ \inprod{\H}{\Delta\X} = 0, \ 
\inprod{\Q}{\Delta\X} \leq 0. 
\end{eqnarray*}
By $O \not=\Delta\X  \in \coneJ$ and 
the assumption that $\Q - \H t \in \mbox{int}\coneJ^*$ for some $t \in \Real$, 
we see that 
\begin{eqnarray*}
0 & < & \inprod{\Q - \H t}{\Delta\X} 
= \inprod{\Q}{\Delta\X} - \inprod{\H}{\Delta\X} t \leq 0,   
\end{eqnarray*}
which is a contradiction. Hence we have shown that both $S(\coneJ)$ and $S(\coneK\cap\coneJ)$ 
are nonempty closed and bounded. 
Therefore, both COP($\coneJ$) and COP($\coneK\cap\coneJ$) have 
optimal solutions, say 
$\widehat{\X}$ and $\X^*$, respectively.  
It follows that $-\infty < \zeta_p(\coneJ) \leq \zeta_p(\coneK\cap\coneJ) = \inprod{\Q}{\X^*} < \infty$. 
Since the equivalence relation~\eqref{eq:equivCond} holds by  Corollary~\ref{corollary:feasiblity}, 
we see that 
$
\zeta_p(\coneJ) =  \zeta_p(\coneK\cap\coneJ) = \inprod{\Q}{\X^*}.   
$ 
Therefore, $\X^*$ is a common optimal solution of COP($\coneJ$) and 
COP($\coneK\cap\coneJ$). 
All other assertions in~\eqref{eq:commonSol} follow from Theorem~\ref{theorem:duality}.
\qed

\section{On Condition (B)}

\label{Section3}

If $\coneJ$ is a face of co$\coneK$ or the convex hull of the union of a family of 
faces of co$\coneK$, then Condition (B) holds. The former case was 
studied throughly in \cite{KIM2020} for its applications to convex conic 
reformulation of QOPs and POPs. In this section, we further investigate 
fundamental properties of Condition (B). 

To characterize Condition (B), we define 
\begin{eqnarray}
\widehat{\FC}(\coneK) & = & \left\{ \coneJ : 
\begin{array}{l}
\coneJ \ \mbox{satisfies Condition (B)},\ {\it i.e.}, \\     
\mbox{$\coneJ$ is a convex cone in co$\coneK$ satisfying co$(\coneK\cap\coneJ) = \coneJ$} 
\end{array}
\right\}. \label{eq:ConditionB} 
\end{eqnarray}

\subsection{Illustrative examples}

We show $4$ examples of $\coneJ \in \widehat{\FC}(\coneK_r)$ for $\coneK_r$ $(r=1,2)$ 
whose convex hull  forms a common semicircular cone in 
the $3$-dimensional space 
in Figure~2, where  a $2$-dim. section of $\coneK_r$ is illustrated $(r=1,2)$. 
We identify a set $S$ on the $2$-dim. section of co$\coneK_r$ with  the 3-dim.
cone$S = \{(\lambda,\lambda \X): \X \in S, \ \lambda \geq 0\}$ (the cone generated by $S$). 
In particular,  an extreme point (or a $1$-dimensional face, respectively) of the section of co$\coneK_r$
corresponds to an extreme ray (or a $2$-dimensional face, respectively) of the semicircular cone co$\coneK_r$.
$\coneK_1$ consists of all extreme rays of the semicircular cone, which correspond to 
the half circle.  
$\coneK_2$ includes the 2-dim. face of 
the semicircular cone, which corresponds to the line segment $[\e,\f]$, 
in addition to all extreme rays. Note that the common $\coneK = \coneK_1$ is used for 
Examples~\ref{example:31}, \ref{example:32} and~\ref{example:33}, and $\coneK = \coneK_2$ for 
Example~\ref{example:34}.
In each example, it is easy to verify that 
assertions (i), (ii) and (iii) of Theorem~\ref{theorem:characterization} with $\coneK = \coneK_r$ 
hold. 

\begin{figure}[t!]  
\begin{center}
\includegraphics[width=0.32\textwidth]{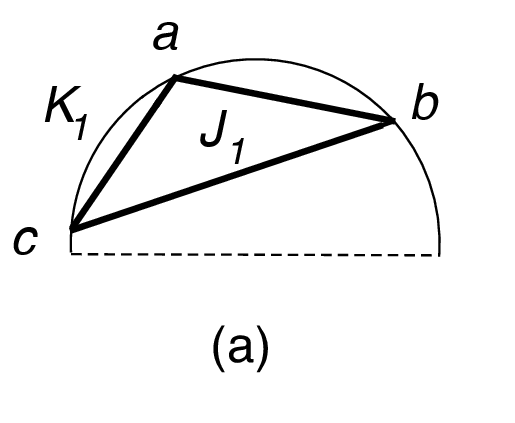}
\mbox{\ } \hspace{2cm}
\includegraphics[width=0.32\textwidth]{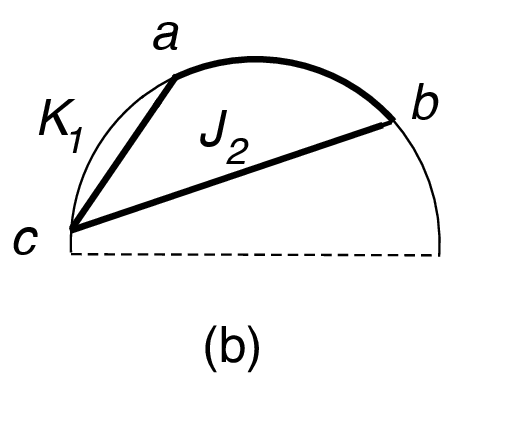}
\includegraphics[width=0.32\textwidth]{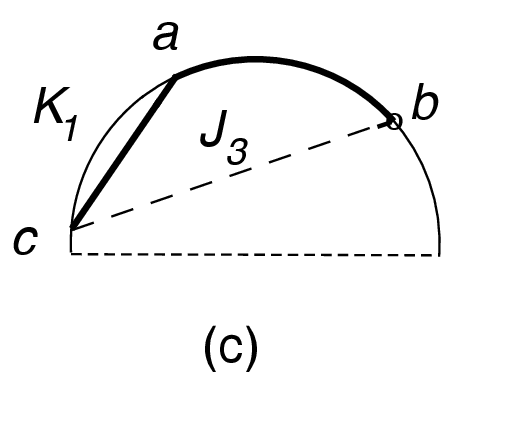}
\mbox{\ } \hspace{2cm}
\includegraphics[width=0.32\textwidth]{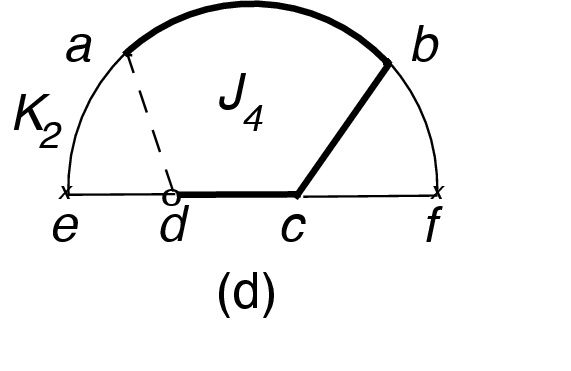}

\end{center}

\caption{
A $2$-dim. section of a $3$-dim. semicircular cone co$\coneK$. 
(a) --- Example 3.1, (b)~---~Example~3.2, (c) --- Example 3.3 and 
(d) --- Example 3.4. 
In (a), (b) and (c), each point on the solid half circle is identified with an 
extreme ray of the $3$-dim. semicircular cone co$\coneK$. In (d), 
the line segment $[\e,\f]$ is identified with  
the $2$-dim. face of the $3$-dim. semicircular cone co$\coneK$.
}
\end{figure} 

\examp \label{example:31}
If we choose $3$ distinct extreme points $\a, \ \b, \ \cc$ on the half circle as 
in Figure~2~(a),
their convex hull $\coneJ_1 \in \widehat{\FC}(\coneK_1)$ forms a closed polyhedral cone. 
Letting 
$\coneC_1 = \mbox{cone}(\{\a, \ \b, \ \cc\})$, we can write $
\coneJ_1 = \mbox{co} (\coneC_1)$. 
\eexamp

\examp \label{example:32}
Let $\coneC_2$ be the cone generated by the union of all extreme points contained in 
the arc $\a$ to $\b$ along the half circle 
and an extreme point $\cc$ (see Figure~2~(b)), their convex hull 
$\coneJ_2 = \mbox{co}(\coneC_2)$ forms a non-polyhedral closed convex cone. 
We see that $\coneJ_2 \in \widehat{\FC}(\coneK_1)$.
\eexamp

\examp \label{example:33}
We modify Example~\ref{example:32} by letting 
$\coneC_3 = \coneC_2\backslash \mbox{cone}(\{\b\})$ as 
in Figure~2~(c). Then $\coneJ_3 = \mbox{co}(\coneC_3) \in \widehat{\FC}(\coneK_1)$. 
In this case, $\coneJ_3$ 
is not closed. Hence, this example shows that $\widehat{\FC}(\coneK_1)$ contains non-closed 
convex cone in general. 
\eexamp

\examp \label{example:34}
In this example, $\coneK_2$ includes the $2$-dim. face of co$\coneK_2$, which 
corresponds to the line segment 
$[\e,\f]$ of its section as in Figure~2~(d). Let $\coneC_4$ be 
the cone generated by the union of 
all extreme points contained in the arc $\a$ to $\b$ along the half circle and 
the semi-closed interval $(\d,\cc]$ on $[\e,\f]$. Then 
$\coneJ_4 = \mbox{co}(\coneC_4) \in \widehat{\FC}(\coneK_2)$. 
It should be noted that 
$\coneJ_4 \in \widehat{\FC}(\coneK_2)$ cannot be generated 
in Examples 3.1, 3.2 and 3.3 since the $2$-dim. face cone$([\e,\f])$ is not included 
in $\coneK_1$.
\eexamp

\subsection{Basic properties on $\widehat{\FC}(\coneK)$}

We now show some fundamental properties on $\widehat{\FC}(\coneK)$ including 
the ones observed in the examples above. 
For a family $\SC$ of subsets of $\spaceV$ 
and a subset $T$ of $\spaceV$, we use notation 
$\bigcup \SC$ to denotes the 
union of all $U \in \SC$, {\it i.e.}, $\bigcup \SC = \bigcup_{U \in \sSC} U$, and 
$T\cap\SC$ to denote $\{T \cap U : U \in \SC\}$. 

\theo \label{theorem:characterization} 
The following assertions hold. 
\begin{description}
\item{(i) } 
$\coneF \in \widehat{\FC}(\coneK)$ 
for every face $\coneF$ of $\coneJ \in \widehat{\FC}(\coneK)$. 
\item{(ii) } $\coneJ \cap \coneF \in \widehat{\FC}(\coneK)$ 
for every face $\coneF$ of co$\coneK$ and $\coneJ \in \widehat{\FC}(\coneK)$. 
\item{(iii) } 
$\coneF \subseteq \coneK$ for every extreme ray $\coneF$ of $\coneJ \in \widehat{\FC}(\coneK)$.
\item{(iv) } $\bigcup(\coneK\cap\FC)
\subseteq \coneK \cap \mbox{co}\big(\bigcup \FC\big) \subseteq \mbox{co}\big(\bigcup \FC\big)$ 
for every $\FC \subseteq \widehat{\FC}(\coneK)$. 
\item{(v) } $\mbox{co}(\bigcup(\coneK\cap\FC)) 
= \mbox{co}(\coneK \cap \mbox{co}\big(\bigcup \FC\big)) = \mbox{co}\big(\bigcup \FC\big)$  
(hence $\mbox{co}\big(\bigcup \FC\big) \in \widehat{\FC}(\coneK)$) 
for every $\FC \subseteq \widehat{\FC}(\coneK)$. 
\end{description} 
\etheo

\noindent
{\em Proof of Theorem~\ref{theorem:characterization} (i): }
Let $\coneF$ be a face of $\coneJ \in \widehat{\FC}(\coneK)$. 
co$(\coneK\cap\coneF) \subseteq \coneF$ is obvious. To show the converse inclusion relation,
let $\X \in \coneF \subseteq \coneJ = \mbox{co}(\coneK \cap \coneJ)$. 
Then $\X = \sum_{k=1}^m\X_k$ 
for some $\X_k \in \coneK \cap \coneJ$ $(1 \leq k \leq m)$. 
Since $\coneF$ is a face of $\coneJ$, 
$\X_k \in \coneF$ $(1 \leq k \leq m)$. 
Hence, $\X \in \mbox{co}(\coneK \cap \coneF)$, and we have 
shown that co$(\coneK\cap\coneF) \supset \coneF$.  
\qed


Obviously co$\coneK \in \widehat{\FC}(\coneK)$. Taking $\coneJ = \mbox{co}\coneK$ in (i), 
we see that $\coneF \in \widehat{\FC}(\coneK)$ for every face $\coneF$ of co$\coneK$. 
In particular, every extreme ray of co$\coneK$ is in $\widehat{\FC}(\coneK)$. 


\noindent
{\em Proof of Theorem~\ref{theorem:characterization} (ii): }
Assume that $\coneF$ is a face of co$\coneK$ and $\coneJ \in \widehat{\FC}(\coneK)$. Then
$\coneJ \cap \coneF$ is a convex cone and 
$\coneK\cap\coneJ\cap \coneF \subseteq \coneJ\cap \coneF$. Hence 
co$(\coneK\cap\coneJ\cap \coneF) \subseteq  
\coneJ\cap \coneF$ follows. 
To show the converse inclusion co$(\coneK\cap\coneJ\cap \coneF) \supset \coneJ\cap \coneF$,  
suppose that $\X \in \coneJ\cap \coneF$. 
Since $\coneJ = \mbox{co}(\coneK\cap\coneJ)$ by assumption, 
$\coneF \ni \X = \sum_{k=1}^m \X_k$ for some 
$\X_k \in \coneK\cap \coneJ \subseteq \mbox{co}\coneK$ $(1 \leq k \leq m)$. 
Since $\coneF$ is a face of co$\coneK$, $\X_k \in \coneF$ $(1 \leq k \leq m)$. 
Hence $\X = \sum_{k=1}^m \X_k \in \mbox{co}(\coneK \cap \coneJ \cap \coneF)$. 
Thus, we have shown co$(\coneK\cap\coneJ\cap \coneF) =   
\coneJ\cap \coneF$ and $\coneJ \cap \coneF \in \widehat{\FC}(\coneK)$.\qed


\noindent
{\em Proof of Theorem~\ref{theorem:characterization} (iii): }
Assume that $\coneF$ is an extreme ray of $\coneJ \in \widehat{\FC}(\coneK)$.
To show $\coneF \subseteq \coneK$, choose an arbitrary nonzero 
$\X \in \coneF$. 
Then, $\X \in \coneF \subseteq \coneJ = \mbox{co}(\coneK \cap \coneJ)$ is represented as 
$\X = \sum_{k=1}^m \X_k$ for some nonzero $\X_k \in \coneK \cap \coneJ$ $(k=1,\ldots,m)$. 
Since $\coneF$ is an extreme ray of $\coneJ$, nonzero $\X_k$ $(1 \leq k \leq m)$ 
all lie in the extreme ray $\coneF$. 
Hence, 
$\X = \lambda_k \X_k$ for some $\lambda_k > 0$ $(1 \leq k \leq m)$. 
Let $k$ be arbitrarily fixed.  Since $\coneK$ is a cone and $\X_k \in \coneK$,  
$\X = \lambda_k \X_k \in \coneK$. Thus, we have shown $\coneF \subseteq \coneK$.
\qed


\noindent
{\em Proof of Theorem~\ref{theorem:characterization} (iv): }
Assume that $\FC \subseteq \widehat{\FC}(\coneK)$. 
Let $\X \in \bigcup (\coneK \cap \FC)$. Then, there exists a $\coneJ \in \FC$ 
such that $\X \in \coneK \cap \coneJ$, which implies that 
$\X \in \coneK \cap \mbox{co}\big(\bigcup\FC\big)$. Hence, 
$\bigcup (\coneK \cap \FC) \subseteq  \coneK \cap \mbox{co}\big(\bigcup\FC\big)$ holds. 
The second inclusion relation is straightforward. 
\qed


\noindent
{\em Proof of Theorem~\ref{theorem:characterization} (v): } 
Assume that $\FC \subseteq \widehat{\FC}(\coneK)$. 
By (iv), it suffices to show that co$(\bigcup (\coneK \cap \FC)) 
\supset \mbox{co}\big(\bigcup \FC\big)$. 
By assumption, 
$\mbox{co}\big(\bigcup (\coneK \cap \FC)\big) \supset \mbox{co}(\coneK \cap \coneJ)=\coneJ $ 
for every $\coneJ \in \FC$. 
Hence, $\mbox{co}\big(\bigcup (\coneK \cap \FC)\big) \supset \bigcup\FC$ follows. 
Since $\mbox{co} \big(\bigcup (\coneK \cap \FC)\big)$ is a convex cone, we obtain that   
$\mbox{co} \big(\bigcup (\coneK \cap \FC)\big) \supset \mbox{co}(\bigcup\FC)$.  
\qed

\subsection{Proof of Theorem~\ref{theorem:NS}} 

\label{Subsection:proofNS}

`only if part' of (i): Assume that $\coneJ \in \widehat{\FC}(\coneK)$. 
Let $\coneK' = \coneK \cap \coneJ$. Then, $\coneK'$ is a cone in $\coneK$ and 
$\coneJ = \mbox{co}\big(\coneK \cap \coneJ\big) = \mbox{co}\coneK'$.

`if part' of (i): Assume that $\coneJ = \mbox{co}\coneK'$ 
for some cone $\coneK' \subseteq \coneK$. 
Since $\coneJ$ is convex, we obviously see that 
$\mbox{co}\big(\coneK \cap \coneJ\big) \subseteq \coneJ$. We also see 
from $\coneK' \subseteq \coneK$ and $\mbox{co}\coneK' = \coneJ$ that 
$\coneK' = \coneK'\cap\mbox{co}\coneK' \subseteq \coneK \cap \coneJ$. Hence 
$\coneJ = \mbox{co}\coneK' \subseteq \mbox{co}\big(\coneK \cap \coneJ\big)$. 
Therefore, we have shown $\mbox{co}\big(\coneK \cap \coneJ\big) = \coneJ$ and 
$\coneJ \in \widehat{\FC}(\coneK)$. 

`only if part' of (ii): Assume that $\coneJ \in \widehat{\FC}(\coneK)$. 
Let $\FC = \{\coneJ\}$. Then, $\FC \subseteq \widehat{\FC}(\coneK)$ and 
$\coneJ = \mbox{co}\coneJ = \mbox{co}\big(\bigcup \FC \big)$ holds since $\coneJ$ is convex. 

`if part' of (ii): Assume that $\coneJ = \mbox{co}\big(\bigcup \FC \big)$ for some 
$\FC \subseteq \widehat{\FC}(\coneK)$. Then, $\coneJ \in \widehat{\FC}(\coneK)$ follows from 
Theorem~\ref{theorem:characterization} (v).

Since $\coneJ \subseteq \mbox{co}\coneK$ and $\H \in \mbox{int}\coneK^*$ by the assumption, 
the feasible region $\{\X \in \coneK\cap\coneJ : \inprod{\H}{\X} = 1 \}$ of 
COP($\coneK\cap\coneJ$), the feasible region 
$\{\X \in \coneJ : \inprod{\H}{\X} = 1 \}$ of COP($\coneJ$) and 
the feasible region 
$S \equiv \{ \X \in \mbox{co}(\coneK\cap\coneJ) : \inprod{\H}{\X} = 1 \}$ 
of COP($\mbox{co}(\coneK\cap\coneJ)$) 
are all bounded and closed.

`only if part' of (iii):
Assume that $\coneJ \in \widehat{\FC}(\coneK)$. Let $\Q$, arbitrarily chosen from $ \spaceV$, be fixed.  
The feasible region 
of COP$(\coneJ)$ is either empty, or 
closed and bounded.
If it is empty, then $\zeta_p(\coneJ) = \zeta_p(\coneK \cap \coneJ) = \infty$. Otherwise,  
it is nonempty, closed and bounded. Hence, 
$-\infty < \zeta_p(\coneJ) <
\infty$. 
Therefore, Conditions (A)' and (B) hold, and 
$\zeta_p(\coneK \cap \coneJ) = \zeta_p(\coneJ)$ follows from 
Corollary~\ref{corollary:feasiblity}. 

`if part' of (iii): 
Assuming $\coneJ \not\in \widehat{\FC}(\coneK)$, we show that 
$-\infty < \zeta_p(\coneJ) < \zeta_p(\coneK \cap \coneJ) < \infty$ 
for some $\Q \in \spaceV$. It follows 
from $\coneJ \not\in \widehat{\FC}(\coneK)$
that the 
convex cone 
$\mbox{co}(\coneK \cap \coneJ)$ is a proper subset of the closed convex cone $\coneJ$. 
Hence, there exists a nonzero 
$\overline{\X}\in\coneJ\backslash\big(\mbox{co}(\coneK\cap\coneJ)\big)\subseteq\mbox{co}\coneK$. 
Let $\widetilde{\X} = \overline{\X}/\inprod{\H}{\overline{\X}}\in\coneJ\backslash\big(\mbox{co}(\coneK\cap\coneJ)\big)\subseteq\mbox{co}\coneK$. Then, $\widetilde{\X}$ is a feasible solution 
of COP$(\coneJ)$ but not in the bounded and closed feasible region 
$S = \{ \X \in \mbox{co}\big(\coneK\cap\coneJ) : \inprod{\H}{\X} = 1 \}$ of COP($\mbox{co}(\coneK\cap\coneJ)$). 
By the separation theorem of 
convex sets (see, for example, \cite[Theorem 11.4.1]{ROCKAFELLAR1970}), there exist 
$\Q \in \spaceV$ such that 
$\inprod{\Q}{\widetilde{\X}} < \inf\{\inprod{\Q}{\X}: \X \in S\} 
= \zeta_p(\mbox{co}(\coneK \cap \coneJ))$. 
Therefore, we obtain that 
\begin{eqnarray*}
-\infty < \zeta_p(\coneJ) \leq \inprod{\Q}{\widetilde{\X}} < 
\zeta_p(\mbox{co}(\coneK \cap \coneJ)) \leq \zeta_p(\coneK \cap \coneJ). 
\end{eqnarray*}
\qed  

Given a convex cone $\coneC$, there are various ways to represent $\coneC$ 
as the convex hull of a nonconvex cone. For example, when the convex cone $\coneC$ is 
closed and pointed, $\coneC$ can be represented as the convex hull of 
the nonconvex cone consisting of the 
extreme rays of $\coneC$ \cite[Theorem 18.5]{ROCKAFELLAR1970}. 
However, any face of $\coneC$ can be added to the nonconvex cone.  
Thus, the representation of 
$\coneC$ in terms of the convex hull of a nonconvex cone $\coneK$ is not unique. 
Theorem~\ref{theorem:NS} (i) shows the possibility of the `finest' 
representation of $\coneJ \in \widehat{\FC}(\coneK)$, and (ii) the possibility of 
various coarse 
representations of $\coneJ \in \widehat{\FC}(\coneK)$. 
A similar observation can be applied to co$\coneK$; two distinct nonconvex cones 
$\coneK_1 \subseteq \spaceV$ and $\coneK_2 \subseteq \spaceV$ induce a common 
convex cone as their convex hull, co$\coneK_1 = \mbox{co}\coneK_2$, such that 
$\widehat{\FC}(\coneK_1) \not= \widehat{\FC}(\coneK_2)$. 
This fact has been observed in Examples~\ref{example:34}.  

\subsection{Decompositions of COP$(\coneK\cap\coneJ)$}

We now focus on 
Theorems~\ref{theorem:characterization} (v) and~\ref{theorem:NS} (ii). 
Let $\coneJ \in \widehat{\FC}(\coneK)$. Assume that $\coneJ$ is decomposed into 
$\coneF \in \FC \subseteq$ for some $\FC \subseteq \widehat{\FC}(\coneK)$ such that 
$\coneJ = \mbox{co}\big(\bigcup \FC\big)$ as in Theorem~\ref{theorem:NS} (ii). 
In general, $\bigcup(\coneK\cap\FC)$ could be a proper subset of 
$\coneK \cap \mbox{co}\big(\bigcup \FC\big)$ by 
Theorem~\ref{theorem:characterization} (iv). 
By Theorem~\ref{theorem:characterization} (v), however, their convex hulls 
coincide with each other, and 
COP(co$(\bigcup(\coneK\cap\FC))$) and 
COP(co$(\coneK \cap \mbox{co}\big(\bigcup \FC\big))$) both induce a common convex relaxation, 
COP($\mbox{co}\big(\bigcup \FC\big)$). We also know that each $\coneF \in \FC$ satisfies 
co$(\coneK \cap \coneF) = \coneF$; hence COP$(\coneF)$ is a convex relaxation of 
COP$(\coneK \cap \coneF)$. 
The following theorem summarizes the relations of the three pairs of 
COPs and their convex conic relaxation, COP($\bigcup(\coneK\cap\FC)$) 
and COP($\mbox{co}\big(\bigcup \FC\big)$), 
COP($\coneK\cap\mbox{co}\big(\bigcup \FC\big)$) and COP($\mbox{co}\big(\bigcup \FC)$), and COP($\coneK\cap\coneF$) and COP($\coneF$) 
($\coneF \in \FC$). In particular, assertion (iii) of the theorem means that 
 COP$(\coneK \cap \mbox{co}\big(\bigcup \FC\big))$ can be decomposed 
into the family of subproblems COP$(\coneK\cap\coneF)$, which 
is reformulated by COP$(\coneF)$, $(\coneF \in \FC)$.

\theo \label{theorem:equivalence} 
Assume that 
$\FC \subseteq \widehat{\FC}(\coneK)$ and that 
$-\infty < \zeta_p(\mbox{co}\big(\bigcup \FC\big)) < \infty$. 
Then, the following assertions hold. 
\begin{description}
\item{(i) } $\zeta_p(\bigcup(\coneK\cap\FC))=
\zeta_p(\coneK\cap\mbox{co}\big(\bigcup \FC\big))
=\zeta_p(\mbox{co}\big(\bigcup \FC\big))$.
\item{(ii) } $\zeta_p(\coneK\cap\coneF)=\zeta_p(\coneF)$ $(\coneF \in \FC)$. 
\item{(iii) } 
$\zeta_p(\coneK\cap\mbox{co}\big(\bigcup \FC\big))
= \inf\{\zeta_p(\coneK\cap\coneF): \coneF\in\FC\} = \inf\{\zeta_p(\coneF): \coneF\in\FC\}$.
\end{description} 
\etheo
\proof{
(i) The pair of $\coneK$ and $\coneJ = \mbox{co}\big(\bigcup \FC\big)$ 
satisfies Condition (A) by Lemma~\ref{lemma:feasiblity}, and 
Condition (B) by Theorem~\ref{theorem:characterization} (v). Hence 
$-\infty<\zeta_p(\coneK\cap\mbox{co}\big(\bigcup \FC\big))=\zeta_p(\mbox{co}\big(\bigcup \FC\big))<\infty$
by Theorem~\ref{theorem:mainKKT}. 
We now prove the identity $\zeta_p(\bigcup(\coneK\cap\FC))=
\zeta_p(\mbox{co}\big(\bigcup \FC\big))$. First, we show that 
$-\infty < \zeta_p(\bigcup(\coneK\cap\FC)) < \infty$. 
Since $\bigcup(\coneK\cap\FC) \subseteq \coneK\cap\mbox{co}\big(\bigcup \FC\big)$ 
by Theorem~\ref{theorem:characterization} (iv), we see 
that $-\infty < \zeta_p(\coneK\cap\mbox{co}\big(\bigcup \FC\big)) \leq \zeta_p(\bigcup(\coneK\cap\FC))$. 
It remains to show that COP($\bigcup(\coneK\cap\FC)$) is feasible. By Condition (A), 
there exists a feasible solution $\overline{\X}$ of 
COP($\coneK \cap \mbox{co}\big(\bigcup \FC\big)$), which 
satisfies $\overline{\X} \in \coneK$, $\overline{\X} \in \mbox{co}(\bigcup\FC)$ 
and $\inprod{\H}{\overline{\X}} = 1$. From $\overline{\X} \in \mbox{co}(\bigcup\FC)$, 
there exist $\X^k \in \bigcup \FC$ $(1 \leq k \leq m)$ such that 
$\overline{\X} = \sum_{k=1}^m\X^k$. Since $1=\inprod{\H}{\overline{\X}}=\sum_{k=1}^m
\inprod{\H}{\X^k}$, $\inprod{\H}{\X^k} > 0$ for some $k$. Let 
$\widehat{\X} = \X^k / \inprod{\H}{\X^k}$. Then $\widehat{\X} \in \coneK \cap \coneF$ 
for some $\coneF \in \FC$ and $\inprod{\H}{\widehat{\X}}=1$, which implies that 
$\widehat{\X}$ is a feasible solution of COP($\bigcup(\coneK\cap\FC)$). 
Hence, we have shown that $-\infty < \zeta_p(\bigcup(\coneK\cap\FC)) < \infty$. 
Now, we observe that $\big(\bigcup(\coneK\cap\FC)\big)\cap\mbox{co}\big(\bigcup \FC\big) = \bigcup(\coneK\cap\FC)$ 
and that 
co$\big(\bigcup(\coneK\cap\FC)\big)\cap\mbox{co}\big(\bigcup \FC\big) = \mbox{co}\big(\bigcup(\coneK\cap\FC)\big) = \mbox{co}\big(\bigcup\FC\big)$ 
by Theorem~\ref{theorem:characterization} (v). Therefore  
$\zeta_p(\bigcup(\coneK\cap\FC)) = \zeta_p(\bigcup(\coneK\cap\FC))\cap\mbox{co}\big(\bigcup \FC\big))
=\zeta_p(\mbox{co}\big(\bigcup \FC\big))$ 
follows from Theorem~\ref{theorem:mainKKT} with replacing $\coneK$ by $\bigcup(\coneK\cap\FC)$ 
and $\coneJ$ by $\mbox{co}\big(\bigcup \FC\big)$.

(ii) Let $\coneF \in \FC$ be  arbitrarily fixed. Then, 
co$(\coneK\cap\coneF) = \coneF$.
By Lemma~\ref{lemma:feasiblity}, if COP($\coneK\cap\coneF$) is infeasible then $\zeta_p(\coneK\cap\coneF) = \zeta_p(\coneF) = \infty$. Otherwise, 
COP($\coneF$) is feasible; hence 
$ \zeta_p(\coneF) < \infty$. 
We also see that 
$-\infty < \zeta_p(\mbox{co}\big(\bigcup \FC\big)) \leq \zeta_p(\coneF)$  
from $\coneF \subseteq \mbox{co}\big(\bigcup \FC\big)$. By applying 
Corollary~\ref{corollary:feasiblity} with replacing 
$\coneJ$ by $\coneF$, we obtain $\zeta_p(\coneK\cap\coneF)=\zeta_p(\coneF)$. 

(iii) By (i) and (ii), it suffices to show $\zeta_p(\bigcup(\coneK\cap\FC)) 
= \inf\{\zeta_p(\coneK\cap\coneF): \coneF\in\FC\}$. 
Since $\coneK \cap \coneF \subseteq \bigcup(\coneK\cap\FC)$ for every $\coneF \in \FC$, we see 
$\zeta_p(\bigcup(\coneK\cap\FC)) \leq \inf\{\zeta_p(\coneK\cap\coneF): \coneF\in\FC\}$. 
To show the converse inequality, let $\overline{\X}$ be an arbitrary feasible solution of 
COP($\bigcup(\coneK\cap\FC)$) with the objective value 
$\bar{\zeta}_p=\inprod{\Q}{\overline{\X}}$. Then  
$\inprod{\H}{\overline{\X}} = 1$ and $\overline{\X} \in \bigcup(\coneK\cap\FC)$, {\it i.e.}, 
$\overline{\X} \in \coneK \cap \overline{\coneF}$ for some $\overline{\coneF} \in \FC$. 
Hence $\overline{\X}$ is a feasible solution of COP($\coneK \cap \overline{\coneF}$). 
Therefore, 
$ 
\inf\{\zeta_p(\coneK\cap\coneF): \coneF\in\FC\} \leq \zeta_p(\coneK\cap\overline{\coneF}) \leq 
\bar{\zeta}_p. 
$\qed
}


\section{A class of quadratically constrained quadratic programs with multiple 
nonconvex constraints}
 
\label{Section4}

\noindent
Throughout this section, we assume that $\coneK = \{\x\x^T: \x \in \Real^n\}$,  where 
$\Real^n$ denotes the $n$-dimensional linear space of column vectors 
$\x=(x_1,\ldots,,x_n)$. Thus, co$\coneK$ forms the positive semidefinite cone 
$\SymMat^n_+$ in the space $\SymMat^n$ of $n\times n$ symmetric 
matrices. 
Let $\coneJ \subseteq \SymMat^n$ be a closed convex cone 
and $\Q, \ \H \in \SymMat^n$. 
We use COP$(\coneK\cap\coneJ,\Q,\H)$ for COP$(\coneK\cap\coneJ)$ 
and $\zeta_p(\coneK\cap\coneJ,\Q,\H)$ for $\zeta_p(\coneK\cap\coneJ)$ to display 
their dependency on $\Q\in\SymMat^n$ and $\H \in \SymMat^n$. Similarly, 
COP$(\coneJ,\Q,\H)$ for COP$(\coneJ)$, $\zeta_p(\coneJ,\Q,\H)$ 
for $\zeta_p(\coneJ)$, DCOP$(\coneJ,\Q,\H)$ for DCOP$(\coneJ)$, and 
$\zeta_d(\coneJ,\Q,\H)$ for $\zeta_d(\coneJ)$. 


COP$(\coneK\cap\coneJ,\Q,\H)$ represents a general (or extended) quadratically constrained 
quadratic program (QCQP) 
\begin{eqnarray*}
\zeta_p(\coneK\cap\coneJ,\Q,\H) & = & \inf\left\{\x^T\Q\x : 
\x \in \Real^n, \ \x\x^T \in \coneJ, \ 
\x^T\H\x = 1 
\right\},  
\end{eqnarray*}
and COP$(\coneJ,\Q,\H)$ its semidefinite programming (SDP) relaxation. 
Recall that Theorem~\ref{theorem:NS} (i) states 
a necessary and sufficient condition 
\begin{eqnarray*}
\mbox{$\coneJ = \mbox{co}\coneK'$ for some cone  
$\coneK' \subseteq \coneK = \cup \{\mbox{all extreme rays of } \SymMat^n_+\}$} 
\end{eqnarray*}
for $\coneJ \in \widehat{\FC}(\coneK)$,   
and Theorem~\ref{theorem:characterization} (iii) describes a necessary 
condition 
\begin{eqnarray*}
\mbox{every extreme ray of $\coneJ$ lies in $\coneK$} 
\end{eqnarray*}
for $\coneJ \in \widehat{\FC}(\coneK)$. See Figure~2 in Section~3. 
These two conditions 
are independent 
from the description of $\coneJ$.  When applying to QCQPs, however, 
$\coneJ$ is usually 
described in terms of inequalities $\inprod{\B_k}{\X}\leq 0$ and/or equalities 
$\inprod{\B_k}{\X} = 0$ with $\B_k \in \SymMat^n$ $(1 \leq k \leq m)$. In this 
section, we focus on the cases where $\coneJ = \{\X \in \SymMat^n_+ : 
\inprod{\B_k}{\X}\leq 0 \ (1 \leq k \leq m)\}$, and present a sufficient 
condition on $\B_k \in \SymMat^n$ $(1 \leq k \leq m)$ for $\coneJ \in \widehat{\FC}(\coneK)$. 
The condition should be sufficient for the above mentioned conditions to hold.


For each $\B \in \SymMat^n$, let 
\begin{eqnarray*}
& & \coneJ_0(\B) = \{ \X \in \SymMat^n_+ : \inprod{\B}{\X} = 0\} , \  
 \coneJ_-(\B) = \{ \X \in \SymMat^n_+ : \inprod{\B}{\X} \leq 0\}.
\end{eqnarray*}
Let $m$ be a nonnegative integer, $\Q, \ \H, \ \B_1,\ldots,\B_m \in \SymMat^n$, 
and $\coneJ_- = \cap_{k=1}^m \coneJ_-(\B_k)$. We note that 
$\coneJ_- = \SymMat^n_+ \in \widehat{\FC}(\coneK)$ if $m=0$.  
We consider COP$(\coneK\cap\coneJ_-,\Q,\H)$ and its convex relaxation COP$(\coneJ_-,\Q,\H)$. 
We can rewrite COP$(\coneK\cap\coneJ_-,\Q,\H)$ as a 
QCQP:  
\begin{eqnarray}
\zeta_p(\coneK\cap\coneJ_-,\Q,\H) & = & 
\inf\left\{ \x^T\Q\x : 
\begin{array}{l}
\x \in \Real^n, \ \x^T\B_k\x \leq 0 \ (1 \leq k \leq m), \\ 
\x^T\H\x = 1 
\end{array}
\right\}.
\label{eq:QCQPmultipleConst}
\end{eqnarray}
COP$(\coneJ_-,\Q,\H)$ serves as an SDP relaxation of the QCQP~\eqref{eq:QCQPmultipleConst}. 
We establish the following result. 
\theo \label{theorem:QCQPmultipleConst} Assume that 
\begin{eqnarray}
\coneJ_0(\B_k) \subseteq \coneJ_- \equiv \cap_{\ell=1}^m \coneJ_-(\B_\ell) \ (1 \leq k \leq m). \label{eq:assumption}
\end{eqnarray} 
(Recall Figure 1 in Section 1). Then,  
\begin{description}
\item{(i) } $\coneJ_- \in \widehat{\FC}(\coneK)$. 
\item{(ii) } Let $\Q\in\SymMat^n$ and $\H\in\SymMat^n$. 
Then 
$-\infty < \zeta_p(\coneJ_-,\Q,\H) = \zeta_p(\coneK\cap\coneJ_-,\Q,\H) < \infty$ iff 
 $-\infty < \zeta_p(\coneJ_-,\Q,\H) < \infty$.
\end{description}
\etheo


We provide some illustrative examples in Section~4.1,
before presenting a proof of the theorem in Section 4.2. 
If $m=0$ or $1$, then condition~\eqref{eq:assumption} obviously holds. 
The following 
proposition presents sufficient conditions for~\eqref{eq:assumption}, which can be 
algebraically verified.  In particular, condition~\eqref{eq:assumption3} is used in 
the examples.

\prop Suppose that $m \geq 2$.  If 
\begin{eqnarray}
-\B_{\ell} - \B_k \tau \in \SymMat^n_+ \ \mbox{for some } \tau \in \Real \ 
(1 \leq k, \ \ell \leq m, \ k\not=\ell), 
 \label{eq:assumption2} 
\end{eqnarray} 
where $\tau$ can depend on both $k$ and $\ell$, or 
\begin{eqnarray}
\hspace{3mm} \; \; \; \inprod{\B_k+\B_\ell}{\X} \leq 0 \ \mbox{for every } \X \in \SymMat^n_+ \ \mbox{or} -(\B_k+\B_\ell)\in \SymMat^n_+ \ (1 \leq k < \ell \leq m), \label{eq:assumption3} 
\end{eqnarray}
then~\eqref{eq:assumption} holds.  
\eprop 
\proof{
Condition~\eqref{eq:assumption2} can be rewritten as 
\begin{eqnarray*}
-\B_\ell \in \SymMat^n_+ + \{\B_k\tau : \tau \in \Real\} \ (1 \leq k, \ \ell \leq m, \ k\not=\ell).  
\end{eqnarray*}
Also,  condition~\eqref{eq:assumption} can be rewritten as 
\begin{eqnarray*}
-\B_\ell \in \coneJ_0(\B_k)^* & = & \big(\SymMat^n_+ \cap \{\X \in \SymMat^n : \inprod{\B_k}{\X} = 0\}\big)^* \\
& = & 
\mbox{cl}\big(\SymMat^n_+ + \{\B_k\tau : \tau \in \Real\}\big) \ (1 \leq k, \ \ell \leq m, \ k\not=\ell),  
\end{eqnarray*}
where cl$C$ denotes the closure of $C \subseteq \SymMat^n$ (see \cite{PATAKI2007} for the second 
identity). Therefore~\eqref{eq:assumption2} implies~\eqref{eq:assumption}. 
Condition~\eqref{eq:assumption3} is a special case of~\eqref{eq:assumption2} with taking 
$\tau = 1$. Hence \eqref{eq:assumption3} implies~\eqref{eq:assumption}. \qed
}

We note that~\eqref{eq:assumption2} is not necessary for condition of~\eqref{eq:assumption}. 
For example, 
take $n=m=2$, 
$\B_1 = {\scriptsize \begin{pmatrix} 1 & 1 \\ 1 & 0 \end{pmatrix}}$ and 
$\B_2 = {\scriptsize \begin{pmatrix} -1 & 0 \\0 & 0 \end{pmatrix}}$. Then 
$-\B_1-\B_2\tau \not\in \SymMat^2_+$ for any $\tau \in \Real$ but 
\begin{eqnarray*}
\coneJ_0(\B_1)  \subseteq  \SymMat^2_+ = \coneJ_-(\B_2), \ 
\coneJ_0(\B_2)  = 
\left\{ \begin{pmatrix} 0 & 0 \\ 0 & X_{22} \end{pmatrix} : X_{22} \geq 0 \right\} 
\subseteq \coneJ_-(\B_1). 
\end{eqnarray*} 
Hence condition~\eqref{eq:assumption} holds. 
If $\tau \leq 0$ in~\eqref{eq:assumption2}, then $\inprod{\B_{\ell}}{\X} \leq 0$ 
for every $\X \in \SymMat^n_+$ satisfying $\inprod{\B_{k}}{\X} \leq 0$; hence 
the constraint $\inprod{\B_\ell}{\X} \leq 0$ (or $\coneJ_-(\B_{\ell})$) 
is redundant.
If we assume that none of the constraints $\inprod{\B_k}{\X} \leq 0$ 
$(1 \leq k \leq m)$ is redundant, then we can take $\tau > 0$ in~\eqref{eq:assumption2}

In Section 4.3, we prove the following result.
\theo \label{theorem:inequality}
Let $\B \in \SymMat^n$. 
\begin{description}
\item{(i) }
$\coneJ_0(\B) \in \widehat{\FC}(\coneK)$ and $\coneJ_-(\B) \in \widehat{\FC}(\coneK)$ are 
equivalent. 
\item{(ii) } $\coneJ_0(\B) \in \widehat{\FC}(\coneK)$. 
\end{description}
\etheo
\noindent
In Section 4.4, we briefly discuss how multiple equality constraints can be added to 
QCQPs~\eqref{eq:QCQPmultipleConst}.

\subsection{Some examples}  

We present six examples. 
\examp 
If $m=1$, then~\eqref{eq:assumption} is satisfied for any $\B_1 \in \SymMat^n$. 
Let $\Q_0, \ \Q_1, \ \Q_2 \in \SymMat^\ell$. Consider the QCQP 
\begin{eqnarray}
\zeta_{\rm QCQP} & = & \inf\left\{\u^T\Q_0\u : \u\in \Real^\ell, \u^T\Q_1\u \leq 1, \ 
\u^T\Q_2\u \leq 1 \right\}. \label{eq:QCQPex1} 
\end{eqnarray}
This form of QCQP was studied in \cite{POLYAK1998,YE2003}. 
They showed that QCQP~\eqref{eq:QCQPex1} can be solved by its SDP relaxation 
under strong duality of its SDP relaxation, which is not assumed here. 
We can transform QCQP~\eqref{eq:QCQPex1} into  
\begin{eqnarray*}
\eta_{\rm QCQP} & = & \inf\left\{\u^T\Q_0\u : 
\begin{array}{l}
\u\in \Real^\ell,\ u_{\ell+1} \in \Real, \\
\u^T\Q_1\u \leq 1, \ 
\u^T\Q_2\u + u_{\ell+1}^2 = 1 
\end{array}
\right\} \\
& = & \inf\left\{\u^T\Q_0\u : 
\begin{array}{l}
\u\in \Real^\ell,\ u_{\ell+1} \in \Real, \\
 \u^T(\Q_1-\Q_2)\u -u_{\ell+1}^2 \leq 0, \ 
\u^T\Q_2\u + u_{\ell+1}^2 = 1 
\end{array}
\right\} \\
& = & \inf \left\{\x^T\Q\x : \x \in \Real^n,\ \x^T\B_1\x \leq 0, \ \x^T\H\x = 1
\right\} \\
& = &  \zeta_p(\coneJ_-,\Q,\H) \ \mbox{(with $m=1$)}, 
\end{eqnarray*}
where $n = \ell+1$, 
\begin{eqnarray*}
\x = \begin{pmatrix} \u \\ u_{l+1} \end{pmatrix}, \ 
\Q = \begin{pmatrix} \Q_0 & \0 \\ \0^T & 0 \end{pmatrix},
\ 
\B_1 = \begin{pmatrix} \Q_1-\Q_2 & \0 \\ \0^T & -1 \end{pmatrix}, \
\H = \begin{pmatrix} \Q_2 & \0 \\ \0^T & 1 \end{pmatrix}.
\end{eqnarray*}
Thus, condition~\eqref{eq:assumption} is satisfied with $m=1$. 
\eexamp

\examp \label{example:equality} Let $\Q, \ \H,  \ \B \in \SymMat^n$. Consider a QCQP 
with two equality constraints. 
\begin{eqnarray*}
\eta_{\rm QCQP} & = & \left\{\x^T\Q\x : \x^T\B\x = 1, \ \x^T\H\x = 1 \right\} \\
& = & \left\{\x^T\Q\x : \x^T(\B-\H)\x = 0, \ \x^T\H\x = 1 \right\} \
 = \ \zeta_p(\coneJ_0(\B_1),\Q,\H),  
\end{eqnarray*}
where $\B_1 = \B - \H$. 
By Theorem~\ref{theorem:inequality} (ii), $\coneJ_0(\B_1) \in \widehat{\FC}(\coneK)$. We can also rewrite the QCQP as  
\begin{eqnarray*}
\eta_{\rm QCQP} & = & \left\{\x^T\Q\x : \x^T(\B_1)\x \leq 0, \ \x^T(\B_2)\x \leq 0, \ \x^T\H\x = 1 \right\}, 
\end{eqnarray*}
where $\B_2 = -(\B-\H)$. 
Then $-(\B_1 + \B_2) = \O \in \SymMat^n_+$; 
hence \eqref{eq:assumption3} holds with $m=2$. This also shows that 
$\coneJ_0(\B_1) = \coneJ_-(\B_1) \cap \coneJ_-(-\B_1) \in \widehat{\FC}(\coneK)$ for 
every $\B_1 \in \SymMat^n$. 
\eexamp

\examp \label{example:STERN1995}
Let $q_k(\u) = \u^T\Q_k\u + 2\b_k^T\u$ be a quadratic function in $\u \in \Real^{\ell}$, 
where $\Q_k \in \SymMat^\ell$, $\b_k\in\Real^\ell$ $(k=0,1)$. 
Consider a QCQP: 
\begin{eqnarray}
\zeta_{\rm QCQP} & = & \inf\left\{ q_0(\u) : \u \in \Real^\ell, -1 \leq q_1(\u) \leq 1 \right\}.
\label{eq:QEQPex2}
\end{eqnarray}
This type of QCQP was studied in \cite{STERN1995} in connection with indefinite trust region subproblems. See also \cite{JEYAKUMAR2014,YE2003}. 
QCQP~\eqref{eq:QEQPex2} can be rewritten as 
\begin{eqnarray*}
\eta_{\rm QCQP} & = & \inf\left\{ \u^T\Q_0\u + 2\b_0^T\u u_{\ell+1} : 
\begin{array}{l}
\u \in \Real^\ell, \ u_{\ell+1} \in \Real, \ u_{\ell+1}^2 = 1, \\
-\u^T\Q_1\u - 2\b_1^T\u u_{\ell+1} - u_{\ell+1}^2 \leq 0, \\
\u^T\Q_1\u + 2\b_1^T\u u_{\ell+1}-  u_{\ell+1}^2 \leq 0
\end{array}
\right\}\\[3pt]
& = & \inf \left\{ \x^T\Q\x : \x^T\B_1\x \leq 0, \ \x^T\B_2\x \leq 0,\ \x^T\H\x = 1 \right\} \\[3pt] 
& = & \zeta_p(\coneJ_-,\Q,\H) \ \mbox{(with $m=2$)}, 
\end{eqnarray*}
where $n = \ell+1$, 
\begin{eqnarray*}
& & \x = \begin{pmatrix} \u \\ u_{l+1} \end{pmatrix} \in \Real^n, \ 
\Q = \begin{pmatrix} \Q_0 & \b_0 \\ \b_0^T & 0  \end{pmatrix} \in \SymMat^n,
\\[3pt]
& & \B_1 = \begin{pmatrix} -\Q_1 & -\b_1 \\ -\b_1^T & -1  \end{pmatrix} \in \SymMat^n, \
\B_2 = \begin{pmatrix} \Q_1 & \b_1 \\ \b_1^T & -1 \end{pmatrix} \in \SymMat^n, \
\H = \begin{pmatrix} \O & \0 & \\  \0^T & 1 \end{pmatrix} \in \SymMat^n.
\end{eqnarray*}
It is easy to verify that 
\begin{eqnarray}
\inprod{\B_1}{\X} + \inprod{\B_2}{\X} = -2X_{nn} \leq 0 \ \mbox{for every } \X \in \SymMat^n_+.\label{eq:B1B2}
\end{eqnarray}
Therefore, condition~\eqref{eq:assumption3} is satisfied with $m=2$.
\eexamp 
 
\examp \label{example:STERN1995Modified}
We add the constraint $\parallel \u \parallel^2/\gamma \geq \gamma$ to QCQP~\eqref{eq:QEQPex2} in 
Example~\ref{example:STERN1995}, where $\gamma > 0$ is a parameter  determined 
later. 
Then, the resulting QCQP can be written as
\begin{eqnarray}
\eta_{\rm QCQP} 
& = & \inf\left\{ q_0(\u) : \u \in \Real^\ell, -1 \leq q_1(\u) \leq 1, \ \parallel \u \parallel^2/\gamma \geq \gamma \right\} \nonumber \\
& = & \inf \left\{ \x^T\Q\x : \x^T\B_k\x \leq 0 \ (k=1,2,3) ,\ \x^T\H\x = 1 \right\}, \label{eq:QCQPmodified} 
\end{eqnarray}
where $n$, $\x$, $\Q$, $\B_1$, $\B_2$ and $\H$ are the same as in Example~\ref{example:STERN1995} and 
$ 
\B_3 = {\scriptsize 
\begin{pmatrix} -\I/\gamma & \0 \\ \0^T & \gamma \end{pmatrix}} 
\in \SymMat^n$. 
In addition to~\eqref{eq:B1B2}, 
\begin{eqnarray}
\left.
\begin{array}{lcl}
\inprod{\B_1}{\X} + \inprod{\B_3}{\X} & =  &
\inprod{{\scriptsize 
\begin{pmatrix}-\Q_1-\I/\gamma & -\b_1 \\ -\b_1^T & -1+\gamma \end{pmatrix}}
}
{\X} \leq 0 \ \mbox{for every } \X \in \SymMat^n_+, \\[3pt] 
\inprod{\B_2}{\X} + \inprod{\B_3}{\X} & = & 
\inprod{{\scriptsize 
\begin{pmatrix} \Q_1-\I/\gamma & \b_1 \\ \b_1^T & -1+\gamma \end{pmatrix}
}}
{\X} \leq 0 \ \mbox{for every } \X \in \SymMat^n_+ 
\end{array}
\right\} \label{eq:condModified}
\end{eqnarray}
hold if we take a sufficiently small $\gamma > 0$.  
Therefore, condition~\eqref{eq:assumption3} is satisfied with $m=3$.
 
Adding the constraint $\parallel \u \parallel^2/\gamma \geq \gamma$ to QCQP~\eqref{eq:QEQPex2} is interpreted as removing the ball 
$\{\u \in \Real^n : \parallel \u \parallel / \gamma < \gamma\}$ 
with the center $\0$,  which lies the interior of the feasible region, from the 
feasible region. The above result implies that 
if the size of the ball is sufficiently small then we can remove the ball 
from the feasible region without destroying $\coneJ_- \in \widehat{\FC}(\coneK)$. 
For example, suppose that $\Q_1=\O$ and 
$\b_1 = {\scriptsize \begin{pmatrix} \0 \\ 1/2 \end{pmatrix}} \in \Real^{\ell}$. 
Then 
QCQP~\eqref{eq:QCQPmodified} turns out to be 
\begin{eqnarray*}
\eta_{\rm QCQP} 
& = & \inf\left\{ q_0(\u) : \u \in \Real^\ell, \ -1 \leq u_\ell \leq 1, \ \parallel \u \parallel^2/\gamma \geq \gamma \right\}.
\end{eqnarray*}
In this case, if $0 < \gamma \leq 4/5$, then~\eqref{eq:condModified} is satisfied. 
It is interesting to note that the unit ball $\{\u \in \Real^\ell : \parallel \u \parallel \leq 1 \}$ is included in $\{\u \in \Real^\ell : -1 \leq u_\ell \leq 1 \}$, 
but we cannot take $\gamma = 1$ to satisfy~\eqref{eq:condModified}.  
\eexamp 

\examp \label{example:n} Let $\B_k$ be a matrix in $\SymMat^n$ whose 
elements satisfy 
\begin{eqnarray*}
[B_k]_{ij} = [B_k]_{ji} \in 
\left\{
\begin{array}{ll}
(-\infty,1] & \mbox{if $i=j=k$}, \\
{(-\infty,-2]}  & \mbox{if $i=j\not=k$}, \\ 
{[-1/(2n),1/(2n)]} & \mbox{otherwise}. 
\end{array}
\right.
\end{eqnarray*}
$(k=1,\ldots,n)$.  
Then, 
\begin{eqnarray*}
-([B_k]_{ij}+[B_\ell]_{ij})= -([B_k]_{ji}+[B_\ell]_{ji}) \in 
\left\{
\begin{array}{ll}
[1,\infty) & \mbox{if $i=j$}, \\
{[-1/n,1/n]} & \mbox{otherwise},  
\end{array}
\right.
\end{eqnarray*}
which implies that $-(\B_k + \B_\ell)$ is diagonally dominant; hence 
positive semidefinite $(1 \leq k < \ell \leq n)$. 
Therefore, condition~\eqref{eq:assumption3} is satisfied with $m=n$.
\eexamp

\examp \label{eq:linearEquality} 
Let $\A$ be an $r \times n$ matrix. Adding a homogeneous linear equality 
constraint $\A\x = \0$ to QCQP~\eqref{eq:QCQPmultipleConst}, we have 
\begin{eqnarray*}
\bar{\eta}_{\rm QCQP} & = & 
\inf\left\{ \x^T\Q\x : 
\begin{array}{l}
\x \in \Real^n, \ \x^T\B_k\x \leq 0 \ (1 \leq k \leq m), \\ 
\x^T\H\x = 1, \ \A\x = \0
\end{array}
\right\} \\
& = & \inf\left\{\inprod{\Q}{\X}:\X \in \coneK\cap\coneJ_-\cap\coneF\right\}
\ = \ \zeta_p(\coneK\cap\coneJ_-\cap\coneF,\Q,\H), 
\end{eqnarray*}
where 
$ 
\coneF \equiv \{ \X \in \SymMat^n_+ : \inprod{\A^T\A}{\X} = 0\}  
$ 
forms a face of $\SymMat^n_+$ since $\A^T\A \in \SymMat^n_+$.  
By Theorem~\ref{theorem:characterization}~(ii), 
$\coneJ_-\cap\coneF \in \widehat{\FC}(\coneK)$ if $\coneJ_- \in \widehat{\FC}(\coneK)$. Therefore, 
we can add $\A\x = \0$ to any of the examples above so that the resulting 
QCQP can still be solved exactly by its SDP relaxation as long as 
$-\infty < \zeta_p(\coneJ_-\cap\coneF,\Q,\H) < \infty$. 
\eexamp

\subsection{Proof of Theorem~\ref{theorem:QCQPmultipleConst}}

We need the following lemma.
\lemm \label{lemma:YeZhang}
(\cite[Lemma 2.2]{YE2003}, see also \cite[Proposition 3]{STURM2003}) Let $\B \in \SymMat^n$ and $\X \in \SymMat^n_+$ with 
rank$\X = r$. Suppose that $\inprod{\B}{\X} \leq 0$. Then, there exists a rank-1 
decomposition of $\X$ such that $\X = \sum_{i=1}^r\x_i\x_i^T$ and $\x_i^T\B\x_i \leq 0$ 
$(1 \leq i \leq r)$. If, in particular, $\inprod{\B}{\X} = 0$, then 
$\x_i^T\B\x_i = 0$ $(1 \leq i \leq r)$. 
\elemm


\noindent
{\em Proof Theorem~\ref{theorem:QCQPmultipleConst} (i)}: 
For the proof, we utilize Theorem~\ref{theorem:NS} (iii) and 
Lemma~\ref{lemma:YeZhang}. 
Choose a $\Q \in \SymMat^n$ arbitrarily, and consider COP$(\coneJ_-,\Q,\I)$. 
We first observe that the feasible region of COP$(\coneJ_-,\Q,\I)$ is either 
empty or bounded since every feasible $\X$ satisfies $\X \in \SymMat^n_+$ 
and $\inprod{\I}{\X} = 1$. 
If $\coneJ_- = \{\O\}$, then the feasible region of COP$(\coneJ_-,\Q,\I)$ is 
empty; hence $\zeta_p(\coneJ_-,\Q,\I) = \zeta_p(\coneK\cap\coneJ_-,\Q,\I) = \infty$. 
Otherwise, there exists a nonzero $\X \in \coneJ_- \subseteq \SymMat^n_+$, 
and $\X/\inprod{\I}{\X}$ lies in the feasible region. Hence, the feasible 
region is bounded, and COP$(\coneJ_-,\Q,\I)$ has a nonzero optimal solution with 
a finite optimal value $\zeta_p(\coneJ_-,\Q,\I)$. 
Obviously $\I \in \SymMat^n_+ 
\subseteq  \coneJ_-^*$. By Theorem~\ref{theorem:duality}, DCOP$(\coneJ_-,\Q,\I)$
has an optimal solution $(\bar{t},\overline{\Y}) \in \Real\times\coneJ_-^*$ such that 
\begin{eqnarray*}
0 & = & \zeta_p(\coneJ_-,\Q,\I) - \zeta_d(\coneJ_-,\Q,\I) = \inprod{\X}{\overline{\Y}} 
\end{eqnarray*}
for every optimal solutions $\X$ of COP$(\coneJ_-,\Q,\I)$.  
The following two cases occur. 
Case (a): there exists a nonzero optimal solution $\X$ 
and a $k \in \{1,\ldots,m\}$ such that $\X \in \coneJ_0(\B_k)$, {\it i.e.}, 
$\inprod{\B_k}{\X} = 0$. Case (b): $\inprod{\B_k}{\X} < 0$ for all optimal solutions 
of COP$(\coneJ_-,\Q,\I)$ and all $k \in \{1,\ldots,m\}$.  

Case (a): 
Let $r = \mbox{rank}\X$. By Lemma~\ref{lemma:YeZhang}, 
there exists a rank-1 decomposition of $\X$ such that 
$\X = \sum_{i=1}^r\x_i\x_i^T$ 
and $\x_i^T\B_k\x_i = 0$ 
({\it i.e.}, $\x_i\x_i^T \in \coneJ_0(\B_k)$) $(1 \leq i \leq r)$. 
By assumption~\eqref{eq:assumption}, $\x_i\x_i^T \in \coneJ_- \ (1 \leq i \leq r)$.
Since 
$1 = \inprod{\I}{\X} = \sum_{i=1}^r \x_i^T\x_i$, there exist a 
$\tau \in (0,1]$ and a $j \in \{1,\ldots,r\}$ such that $\x_j^T\x_j = \tau$. 
Let $\overline{\X} = \x_j\x_j^T/\tau$. Since $ \x_j\x_j^T\in \coneJ_-$, 
$\overline{\X} \in \coneJ_-$. We also see that $\inprod{\I}{\overline{\X}} = \inprod{\I}{\x_j\x_j^T/\tau} = 1$. Hence $\overline{\X}$ 
is a rank-$1$ feasible solution of COP$(\coneJ_-,\Q,\I)$. Furthermore, 
we see from $\overline{\Y} \in \coneJ_-^*$ and $\x_i\x_i^T \in \coneJ_-$ 
$(1 \leq i \leq r)$ that 
\begin{eqnarray*}
0 \leq \inprod{\overline{\Y}}{\overline{\X}} = \frac{\inprod{\overline{\Y}}{\x_j\x_j^T}}{\tau} \leq 
\frac{\inprod{\overline{\Y}}{\X}}{\tau} = 0.  
\end{eqnarray*}
Hence, $\overline{\X}$ is a rank-1 optimal solution of COP$(\coneJ_-,\Q,\I)$, 
and it is an optimal solution of COP$(\coneK\cap\coneJ_-,\Q,\I)$ with the 
same objective value $\inprod{\Q}{\overline{\X}}$. 
Therefore, we have shown that 
$\zeta_p(\coneJ_-,\Q,\I) = \zeta_p(\coneK\cap\coneJ_-,\Q,\I)$. 

Case (b): Let $S$ denote the optimal solution set of COP$(\coneJ_-,\Q,\I)$. 
By the assumption of this case, $\inprod{\B_k}{\X} < 0$ $(1 \leq k \leq m)$ 
for all $\X \in S$. Hence $S$ coincides with the optimal solution set of 
COP$(\SymMat^n_+,\Q,\I)$, an SDP of minimizing $\inprod{\Q}{\X}$ 
subject to $\X \in \SymMat^n_+$ and the single equality constraint 
$\inprod{\H}{\X} = 1$. It is well-known that every solvable SDP with $r$ equality 
constraints has an optimal solution $\X$ with rank $\leq \sqrt{2r}$ 
(see, for example, \cite{PATAKI1998}). Hence, there exists an 
$\overline{\X} \in S$ such that rank$\overline{\X} = 1$, which is a common 
optimal solution of COP$(\coneJ,\Q,\I)$ and COP$(\coneK\cap\coneJ_-,\Q,\I)$. 
Therefore, $\zeta_p(\coneJ_-,\Q,\I) = \zeta_p(\coneK\cap\coneJ_-,\Q,\I)$. 

In both Cases (a) and (b), we have shown that 
$\zeta_p(\coneJ_-,\Q,\I) = \zeta_p(\coneK\cap\coneJ_-,\Q,\I)$. Since $\Q$ is chosen arbitrarily from 
$\SymMat^n$,  
we can conclude by Theorem~\ref{theorem:NS} (iii) that 
$\coneJ_- \in \widehat{\FC}(\coneK)$. 
\qed


\noindent
{\em Proof Theorem~\ref{theorem:QCQPmultipleConst} (ii)}: The desired result follows from Corollary~\ref{corollary:feasiblity}. \qed

\subsection{Proof of Theorem~\ref{theorem:inequality}}

{\em (i) }
We first show that $\coneJ_0(\B) \in \widehat{\FC}(\coneK)$ $\Rightarrow$ $\coneJ_-(\B) \in \widehat{\FC}(\coneK)$. 
Since $\coneJ_-(\B)$ is a convex set containing 
$\coneK\cap\coneJ_-(\B)$, co$(\coneK\cap\coneJ_-(\B)) \subseteq \coneJ_-(\B)$ is obvious. To show 
$\coneJ_-(\B) \subseteq \mbox{co}(\coneK\cap\coneJ_-(\B))$,  let $\X \in \coneJ_-(\B)$. Then 
there exists an $\X^i \in \coneK$ such that 
\begin{eqnarray*}
\X = \sum_{i=1}^m \X^i, \ 0 \geq \inprod{\B}{\X} = \sum_{i=1}^m \inprod{\B}{\X^i}.  
\end{eqnarray*}
Let 
$ 
I_+ = \{ i : \inprod{\B}{\X^i} > 0 \} \ \mbox{and } \ I_0 = \{ i : \inprod{\B}{\X^i} = 0 \}, \ 
I_- = \{ i : \inprod{\B}{\X^i} < 0 \}. 
$ 
By definition, 
$ 
\X^i \in \coneK \cap \coneJ_-(\B) \ (i \in I_0 \cup I_-)
$. 
Hence, if $I_+ = \emptyset$ then $\X \in \mbox{co}(\coneK\cap\coneJ_-(\B))$. 
Now assume that $I_+ \not= \emptyset$. Then, $I_- \not= \emptyset$ since otherwise 
$\inprod{\B}{\X} = \sum_{i=1}^m \inprod{\B}{\X^i} > 0$, a contradiction. Let 
\begin{eqnarray*}
\X^+ = \sum_{i\in I_+} \X^i \in \mbox{co} \coneK,\  \
\X^0 = \left\{\begin{array}{ll} \displaystyle
\sum_{i\in I_0} \X^i & \mbox{if } I_0 \not= \emptyset, \\ 
\0 & \mbox{otherwise }
\end{array}, \right. \ \
\X^- = \sum_{i\in I_-} \X^i.
\end{eqnarray*}
Then, 
\begin{eqnarray}
& & 
\X^0 \in \mbox{co}(\coneK\cap\coneJ_0(\B)), \  \X^- \in \mbox{co}(\coneK\cap\coneJ_-(\B)), \ 
\X = \X^+ + \X^0 + \X^-, \label{eq:3points} \\
& & \alpha^+ \equiv \inprod{\B}{\X^+} > 0, \ \alpha^0 \equiv \inprod{\B}{\X^0} = 0, \ 
\alpha^- \equiv - \inprod{\B}{\X^-} > 0, \nonumber \\  
& & 0 \geq \inprod{\B}{\X} = \inprod{\B}{\X^+ + \X^0 + \X^-} 
= \alpha^+ - \alpha^-. 
\label{eq:alpha}
\end{eqnarray}
Define 
\begin{eqnarray}
\overline{\X} & \equiv & \frac{\alpha^-\X^+ + \alpha^+\X^-}{\alpha^++\alpha^-}, 
\ \mbox{which implies } 
\X^+  = \frac{\alpha^++\alpha^-}{\alpha^-}\overline{\X}- \frac{\alpha^+}{\alpha^-}\X^-. 
\label{eq:Xplus} 
\end{eqnarray}
Then, 
\begin{eqnarray*}
& & \overline{\X} \in \mbox{co}\coneK,\   
\inprod{\B}{\overline{\X}} = 
\frac{\alpha^-\inprod{\B}{\X^+} + \alpha^+\inprod{\B}{\X^-}}{\alpha^++\alpha^-} 
= \frac{\alpha^-\alpha^+ - \alpha^+\alpha^-}{\alpha^++\alpha^-} = 0. 
\end{eqnarray*}
Hence, 
$ 
\overline{\X} \in 
\coneJ_0(\B) = \mbox{co}(\coneK \cap \coneJ_0(\B)) 
\subseteq \mbox{co}(\coneK \cap \coneJ_-(\B)) 
$.  
It follows from \eqref{eq:3points} and \eqref{eq:Xplus} that 
\begin{eqnarray*}
\X  =  \X^+ + \X^0 + \X^- 
 = \frac{\alpha^++\alpha^-}{\alpha^-}\overline{\X} + 
 \frac{\alpha^- - \alpha^+}{\alpha^-}\X^- + \X^0. 
\end{eqnarray*}
Since we have already shown that the three points $\overline{\X}, \X^0$ and $\X^-$ 
lies in co$(\coneK\cap\coneJ_-(\B))$ and $\frac{\alpha^++\alpha^-}{\alpha^-} \geq 0$ and 
$ \frac{\alpha^- - \alpha^+}{\alpha^-} \geq 0$ (by~\eqref{eq:alpha}), we obtain $\X \in co(\coneK\cap\coneJ_-(\B))$. 
Thus we have shown that $\coneJ_0(\B) \in \widehat{\FC}(\coneK)$ 
$\Rightarrow$ $\coneJ_-(\B) \in \widehat{\FC}(\coneK)$. 

We show that   $\coneJ_-(\B) \in \widehat{\FC}(\coneK)$ $\Rightarrow$ $\coneJ_0(\B) \in \widehat{\FC}(\coneK)$. 
It is easy to verify that $\coneJ_0(\B)$ is a face of $\coneJ_-(\B) \in \widehat{\FC}(\coneK)$. 
Then $\coneJ_0(\B) \in \widehat{\FC}(\coneK)$ follows by Theorem~\ref{theorem:characterization} (i). Therefore, we have shown that $\coneJ_0(\B) \in \widehat{\FC}(\coneK)$ 
$\Leftrightarrow$ $\coneJ_-(\B) \in \widehat{\FC}(\coneK)$. 
\qed


\noindent
{\em (ii) } Since $\coneJ_0(\B) \in \widehat{\FC}(\coneK)$ is equivalent to $\coneJ_-(\B) 
\in \widehat{\FC}(\coneK)$ by (i), $\coneJ_0(\B) \in \widehat{\FC}(\coneK)$ follows from 
Theorem~\ref{theorem:QCQPmultipleConst}~(i) with $m=1$. \qed

\rema Assertion (i) of Theorem~\ref{theorem:inequality} holds for a more general case 
where $\coneK$ is a cone in a finite dimensional space $\spaceV$ 
and $\B \in \spaceV$. Define 
\begin{eqnarray*}
\coneJ_0(\B) = \{\X \in \mbox{co}\coneK : \inprod{\B}{\X} = 0\}, \ 
\coneJ_-(\B) = \{\X \in \mbox{co}\coneK : \inprod{\B}{\X} \leq 0\}.
\end{eqnarray*}
Then, $\coneJ_0(\B) \in \widehat{\FC}(\coneK)$ and 
$\coneJ_-(\B) \in \widehat{\FC}(\coneK)$ are equivalent. The proof of 
Theorem~\ref{theorem:inequality} (i) stated above remains valid for this general 
case without any change. 
\erema

\subsection{Adding equality constraints}

Assume 
 that, in general, 
\begin{eqnarray}
\coneJ \equiv \left\{ \X \in \SymMat^n : \inprod{\B_j}{\X}=0 \ (j \in I_0), \ \inprod{\B_k}{\X}\leq0 \ (k \in I_-) \right\} \in \widehat{\FC}(\coneK), 
\label{eq:generalCase} 
\end{eqnarray}
where $I_0$ and $I_-$ denote a partition of $\{1,\ldots,m\}$ and $\Q, \ \H, \ \B_k \in \SymMat^n$ 
$(k \in I_0\cup I_-)$.  
We present a method for adding an equality constraint $\inprod{\B_{m+1}}{\X} = 0$  
to $\coneJ$ so that $\coneJ \cap \coneJ_0(\B_{m+1})$ remains in $\widehat{\FC}(\coneK)$. 
The method can be used recursively by replacing $I_0$ with $I_0 \cup \{m+1\}$. 
Note that $I_-$ is not updated. 
Initially, we take $\coneJ_- \equiv 
\cap_{k=1}^m \coneJ_-(\B_k)$ with $\B_k \in \SymMat^n$ $(1 \leq k \leq m)$ satisfying~\eqref{eq:assumption} or $\coneJ_0(\B_1)$ with $\B_1 \in \SymMat^n$ for $\coneJ$.
  
If $I_0 = \emptyset$, we can choose 
any $\B \in \SymMat^n$ satisfying 
$\coneJ_0(\B) \subseteq \coneJ_-(\B_k)$ $(k \in I_-)$ so that $\coneJ_- \cap \coneJ_0(\B) = \coneJ_0(\B) \in \widehat{\FC}(\coneK)$ by Theorem~\ref{theorem:inequality} (i). As a result, all inequality constraints $\inprod{\B_k}{\X}$ $(k \in I_-)$  become 
redundant. We exclude this trivial choice in the subsequent discussion. 

Let $\coneK' = \coneK \cap \coneJ$, which is a cone in $\SymMat^n$.  
It follows from $\coneJ \in \widehat{\FC}(\coneK)$ that 
$\mbox{co}(\coneK'\cap\coneJ) = \mbox{co}(\coneK\cap\coneJ) =   \coneJ$. Hence $\coneJ \in \widehat{\FC}(\coneK')$. 
All the results in Section 3 is quite general so that 
the results can be applied even if 
$\coneK$ is replaced  with $\coneK'$ and $\SymMat^n_+ = \mbox{co}\coneK$ with co$\coneK'$. However, 
it is difficult to extend Theorems~\ref{theorem:QCQPmultipleConst} and~\ref{theorem:inequality} 
to the general case where $\SymMat^n_+$ is replaced with $\coneJ=\mbox{co}\coneK'$ defined by~\eqref{eq:generalCase}. The main reason is that Lemma~\ref{lemma:YeZhang} is no longer 
valid for the general case. Consequently, it becomes clear that  
adding an equality 
constraint to $\coneJ$ is more restrictive than adding inequalities 
as  only 
the general 
results in Section~3 can be used. 
  

In general, 
$\coneJ_0(\B)$ is a face of 
$\coneJ$ iff 
\begin{eqnarray*}
\B \in \coneJ^* = \mbox{cl} \left\{ 
\Y + \sum_{j \in I_- \cup I_0} \B_jy_j  : 
\Y \in \SymMat^n_+, \ y_j \in \Real \ (j \in I_0),  \ y_k \leq 0 \ (k \in I_-) 
\right\}. 
\end{eqnarray*}
(See \cite{PATAKI2007} for the identity). 
Obviously, every $\B \in \SymMat^n_+$ 
lies in $\coneJ^*$. In this case, $\coneJ_0(\B)$ forms a common face of $\coneJ$ 
and $\SymMat^n_+$. It is also straightforward to see that $\B = -\B_\ell \in \coneJ^*$ 
for every $\ell \in I_-$. In this case, 
$\coneJ \cap \coneJ_0(-\B_\ell) = \big(\cap_{j \in I_0} \coneJ_0(\B_j)\big) \cap \coneJ_0(\B_\ell)$ 
and all inequalities $\inprod{\B_k}{\X} \leq 0$ $(k \in I_-)$ become redundant. 
(Recall that assumption~\eqref{eq:assumption} holds 
if $I_- \not= \emptyset$). More generally, if $y_k \leq 0$ $(k \in I_-)$ and  $y_j \in \Real$ $(j \in I_0)$, then then $\B = \sum_{j \in I_- \cup I_0} \B_jy_j \in \coneJ^*$.  

\section{Concluding discussion}

By extending the condition co$(\coneK\cap\coneJ) = \coneJ$ with a face $\coneJ$  
of co$\coneK$ in \cite{KIM2020} to characterizations of 
the family $\widehat{\FC}(\coneK)$ of all convex cone 
 $\coneJ \subseteq \mbox{co}\coneK$ satisfying 
co$(\coneK\cap\coneJ) = \coneJ$, 
we have established the fundamental properties of $\widehat{\FC}(\coneK)$ in Section 3.
In particular, by applying the properties to nonconvex QCQPs, 
we have shown that a new class of QCQP with multiple nonconvex  inequality and equality  constraints can be solved exactly by
its SDP relaxation in Section~4.

The important and distinctive feature of the geometric nonconvex COP$(\coneK \cap \coneJ)$ 
and its convex conic reformulation COP$(\coneJ)$ is independence from the description of 
$\coneK$ and $\coneJ$. 
The required main assumption is that $\coneK$ is a cone in $\spaceV$ and 
$\coneJ \in \widehat{\FC}(\coneK)$, 
which indicates that the results presented 
in Sections~2 and~3 can be applied in various cases. 
It should be noted that the other assumption $-\infty < \zeta_p(\coneJ) < \infty$ 
is necessary and sufficient to ensure 
$\zeta_p(\coneJ) =  \zeta_p(\coneK\cap\coneJ)$ 
under $\coneJ \in \widehat{\FC}(\coneK)$. We have not imposed any assumption on 
the objective function $\inprod{\Q}{\X}$. 
See Corollary~\ref{corollary:feasiblity}.

To the question of whether Theorem~\ref{theorem:QCQPmultipleConst} can be 
applied to the completely positive relaxation of QCQPs,
we should note that the answer is negative,  
as shown in the following simple example:
Let $\coneK = \{\x\x^T : \x \in \Real^2_+\}$ 
and $\coneJ_0 = \{ \X \in \mbox{co}\coneK : X_{11} - X_{22} = 0 \}$. Then 
$\mbox{co}(\coneK \cap \coneJ_0) = 
\left\{\lambda {\scriptsize \begin{pmatrix} 1 & 1 \\ 1 & 1 \end{pmatrix}} :\lambda \geq 0 \right\}$ 
is a proper subset of $\coneJ_0$; hence $\coneJ_0 \not\in \widehat{\FC}(\coneK)$. 
Thus assertion (ii) of Theorem~\ref{theorem:inequality} does not hold for this case. 
Or, at least, some modification is necessary on assumption~\eqref{eq:assumption}.
This example also shows that Lemma~\ref{lemma:YeZhang} 
is no longer valid if  $\SymMat^n_+$ is replaced with the completely positive cone. 

We recently became aware of the results in
 the paper  \cite{ARGUE2023}  by Argue et al., which investigated the equivalence between QCQPs and their SDP relaxations. 
Some results there including their Proposition 1 are 
related to our main result, Theorem~\ref{theorem:QCQPmultipleConst}. 
We should mention that  our approach is different from theirs and our results were obtained independently from theirs. Our approach,  originated from our previous paper \cite{KIM2020} 
which discussed the equivalent 
convex relaxation of a geometric conic optimization problem in a finite dimensional space, is 
geometric in nature. In contrast,  theirs is limited to QCQP in $\Real^n$ and its 
SDP relaxation, and  does not involve geometric considerations.
Their paper \cite{ARGUE2023} was also written independently 
from our earlier paper \cite{KIM2020}. 
Both approaches should be useful for  further development of the subject. 


\end{document}